\documentclass[14pt]{article}
\usepackage{mathrsfs}
\usepackage{amssymb}
\usepackage{amsthm}
\usepackage{amsmath}
\usepackage{graphicx}
\usepackage{color}
\usepackage{amsfonts}
\usepackage{float}
\usepackage{cite}
\usepackage[text={140mm,210mm},left=45mm,vmarginratio=1:1]{geometry}
\newtheorem{theorem}{Theorem}[section]

\newtheorem{lemma}[theorem]{Lemma}
\newtheorem{proposition}[theorem]{Proposition}

\numberwithin{equation}{section}
\normalsize

\begin{document}
\title{\textbf{Survival probabilities of high-dimensional stochastic SIS and SIR models with random edge weights}}

\author{Xiaofeng Xue \thanks{\textbf{E-mail}: xfxue@bjtu.edu.cn \textbf{Address}: School of Science, Beijing Jiaotong University, Beijing 100044, China.}\\ Beijing Jiaotong University}

\date{}
\maketitle

\noindent {\bf Abstract:} In this paper, we are concerned with the stochastic SIS (susceptible-infected-susceptible) and SIR (susceptible-infected-recovered) models on high-dimensional lattices with random edge weights, where a susceptible vertex is infected by an infectious neighbor at rate proportional to the weight on the edge connecting them. All the edge weights are assumed to be i.i.d.. Our main result gives mean field limits for survival probabilities of the two models as the dimension grows to infinity, which extends the main conclusion given in \cite{Xue2017} for classic stochastic SIS model.


\section{Introduction}\label{section 1}
In this paper, we are concerned with the stochastic SIS (susceptible-infected-susceptible) and SIR (susceptible-infected-recovered) models on high-dimensional lattices $\mathbb{Z}^d$. For later use, we introduce some notations. We use $O$ to denote the origin of $\mathbb{Z}^d$. For each $x\in \mathbb{Z}^d$, we denote by $\|x\|$ the $l_1$ norm of $x$, i.e.,
\[
\|x\|=\sum_{i=1}^d|x_i|
\]
for $x=(x_1,\ldots,x_d)$. For $1\leq i\leq d$, we use $e_i$ to denote the $i$th basic unit-vector of $\mathbb{Z}^d$, i.e.,
\[
e_i=(0,\ldots,0,\mathop 1\limits_{i \text{th}},0\ldots,0).
\]
For $x,y\in \mathbb{Z}^d$, we write $x\sim y$ when and only when $\|x-y\|=1$. We use $E_d$ to denote $\Big\{\{x,y\}:~x\sim y\Big\}$, which is identified with the set of edges on $\mathbb{Z}^d$. For any set $A$, we denote by $|A|$ the cardinality of $A$.

Let $\rho$ be a random variable that $P(0\leq\rho\leq \Theta)=1$ for some $\Theta \in (0,+\infty)$ and $P(\rho>0)>0$, then we assume that $\{\rho(e)\}_{e\in E_d}$ are i.i.d. copies of $\rho$. For $e=\{x,y\}\in E_d$, we write $\rho(e)$ as $\rho(x,y)$. Note that $\rho(x,y)=\rho(y,x)$.

When $\{\rho(e)\}_{e\in E_d}$ are given, the stochastic SIS model  with edge weights $\{\rho(e)\}_{e\in E_d}$ is a continuous-time Markov process $\{C_t\}_{t\geq 0}$ with state space
\[
X_1=\{C:~C\subseteq \mathbb{Z}^d\}
\]
and transition rate function given by
\begin{equation}\label{equ 1.1 transition rate for SIS}
C_t\rightarrow
\begin{cases}
C_t\setminus\{x\} \text{~at rate~} 1 &\text{~if~}x\in C_t,\\
C_t\cup\{x\} \text{~at rate~} \frac{\lambda}{2d}\sum\limits_{y:y\sim x}\rho(x,y)1_{\{y\in C_t\}} &\text{~if~}x\not \in C_t,
\end{cases}
\end{equation}
where $\lambda$ is a positive constant called the infection rate and $1_A$ is the indicator function of the event $A$.

The stochastic SIR model with edge weights $\{\rho(e)\}_{e\in E_d}$ is a continuous-time Markov process $\{(S_t,I_t)\}_{t\geq 0}$ with state space
\[
X_2=\{(S,I):~S\subseteq \mathbb{Z}^d, I\subseteq \mathbb{Z}^d, S\cap I=\emptyset\}
\]
and transition rate function given by
\begin{align}\label{equ 1.2 transtion for SIR}
&(S_t,I_t)\rightarrow\\
&\begin{cases}
(S_t,I_t\setminus\{x\}) \text{~at rate~}1 &\text{~if~}x\in I_t,\\
(S_t\setminus\{x\},I_t\cup\{x\}) \text{~at rate~} \frac{\lambda}{2d}\sum\limits_{y:y\sim x}\rho(x,y)1_{\{y\in I_t\}} &\text{~if~}x\in S_t.
\end{cases}\notag
\end{align}

Both the SIS model and the SIR model describe the spread of epidemics on a graph. For the SIS model, each vertex is in one of two states, `susceptible' or `infectious'. $C_t$ is the set of infectious vertices at moment $t$. An infectious vertex waits for an exponential time with rate one to become susceptible while a susceptible vertex is infected by an infectious neighbor at rate proportional to the weight on the edge connecting them. For the SIR model, each vertex is in one of three states, `susceptible', `infectious' or `recovered'. $S_t$ is the set of susceptible vertices and $I_t$ is the set of infectious vertices at the moment $t$. A susceptible vertex is infected in the same way as that of the SIS model while an infectious vertex waits for an exponential time with rate one to become recovered. A recovered vertex can never infect neighbors or be infected again.

The SIS model is also named as the contact process. The classic contact process is introduced by Harris in \cite{Har1974}, where $\rho=1$. For a detailed survey of the classic contact process, see Chapter 6 of \cite{Lig1985} and Part one of \cite{Lig1999}. The contact process with i.i.d edge weights is introduced by Chen and Yao in \cite{Yao2012}, where the complete convergence theorem of the process is proved. When $P(\rho=1)=p=1-P(\rho=0)$, the model reduces to the contact process on clusters of bond percolation, which is also introduced by Chen and Yao in \cite{Chen2009} to prove a similar complete convergence theorem. It is also interesting to put the random weights on vertices instead of edges, where a susceptible vertex $x$ with weight $\rho(x)$ is infected by an infectious neighbor $y$ with weight $\rho(y)$ at rate proportional to $\rho(x)\rho(y)$. This model is introduced by Peterson on the complete graph in \cite{Pet2011}, where a phase transition phenomenon consistent with a mean-field analysis is shown. Xue studies the contact process with random vertex weights on the oriented lattice in \cite{Xue2015}, where a limit theorem of the critical infection rate is given. When the vertex weight takes $1$ with probability $p$ and takes $0$ otherwise, the process reduces to that on clusters of site percolation, which is a special case of the model introduced in \cite{Ber2011} with $N=1$. In \cite{Ber2011}, Bertacchi, Lanchier and Zucca study the contact process on $C_\infty\times K_N$, where $C_\infty$ is unique infinite open cluster of the site percolation on $\mathbb{Z}^d$ while $K_N$ is the complete graph with $N$ vertices. Criteria to judge whether the process survives is given in \cite{Ber2011}.

The initial motivation of the study in this paper is to extend the main result in \cite{Xue2017}, which gives the mean field limit for survival probability of high-dimensional classic contact process, to the case where the contact process is with random edge weights. We find out that the SIR model is a useful auxiliary tool for us to accomplish our objective and similar conclusion holds for the SIR model simultaneously according to our proof. We are inspired a lot by the technique introduced in \cite{Xue2017b}, which gives asymptotic behavior of the critical value of the high-dimensional SIR model on clusters of bond percolation.

\section{Main results}\label{section 2}
In this section we give our main results. First we introduce some notations and definitions. We assume that the edge weights $\{\rho(e)\}_{e\in E_d}$ are defined under the probability space $(\Omega_d,\mathcal{F}_d,\mu_d)$. The expectation operator with respect to $\mu_d$ is denoted by $E_{\mu_d}$. For $\omega\in \Omega_d$, we write $\rho(e)$ as $\rho(e,\omega)$ when we emphasize that the weight on $e$ is with respect to the random environment $\omega$. For $\lambda>0$ and $\omega\in \Omega_d$, we denote by $P_{\lambda,\omega}$ the probability measure of the SIS and SIR models on $\mathbb{Z}^d$ with infection rate $\lambda$ and edge weights $\{\rho(e,\omega)\}_{e\in E_d}$. $P_{\lambda,\omega}$ is called the quenched measure. We define
\[
P_{\lambda,d}(\cdot)=E_{\mu_d}\Big(P_{\lambda,\omega}(\cdot)\Big)=\int P_{\lambda,\omega}(\cdot)~\mu_d(d \omega).
\]
$P_{\lambda,d}$ is called the annealed measure. When we do not need to distinguish the dimension $d$, we omit the subscript $d$ in the above notations. For $A\subseteq \mathbb{Z}^d$, we write
$C_t$ as $C_t^A$ when $C_0=A$. If $A=\{x\}$ for $x\in \mathbb{Z}^d$, then we write $C_t^{\{x\}}$ as $C_t^x$ for simplicity. For any $x\in \mathbb{Z}^d$, we write $(S_t,I_t)$ as $(S_t^x,I_t^x)$ when $(S_0,I_0)=(\mathbb{Z}^d\setminus\{x\},\{x\})$.

The following theorem is our main result, which gives mean field limits of the survival probabilities of the SIS and SIR models as the dimension $d$ grows to infinity.
\begin{theorem}\label{theorem 2.1 main}
Let $O$ be the origin of $\mathbb{Z}^d$ as we have defined in Section \ref{section 1}, then
\[
\lim_{d\rightarrow+\infty}P_{\lambda,d}\big(C_t^O\neq \emptyset,\forall~t\geq 0\big)=\lim_{d\rightarrow+\infty}P_{\lambda,d}\big(I_t^O\neq \emptyset,\forall~t\geq 0\big)=\frac{\lambda E\rho-1}{\lambda E\rho}
\]
for any $\lambda\geq \frac{1}{E\rho}$, where $E\rho$ is the expectation of $\rho$.
\end{theorem}

Theorem \ref{theorem 2.1 main} shows that for high-dimensional SIS and SIR models with random edge weights, assuming that $O$ is the unique infectious vertex at $t=0$ while other vertices are susceptible, then the probability that infectious vertices will never die out approximately equals $(\lambda E\rho-1)/\lambda E\rho$. This result can be intuitively explained according to a mean-field analysis. When the dimension $d$ is large, it is not likely that infectious vertices will cluster, then $|C_t|$ decreases by one at rate $|C_t|$ and increases by one at rate approximate to
\[
\sum_{x\in C_t}\sum_{y:y\sim x}\frac{\lambda}{2d}\rho(x,y)\approx \lambda |C_t| E\rho
\]
according to the law of large numbers. Then, the embedded chain of $|C_t|$ is similar with a biased random walk on $\mathbb{Z}^1$ that increases by one with probability $\frac{\lambda E\rho}{\lambda E\rho+1}$ or decreases by one with probability $\frac{1}{\lambda E\rho +1}$. Such a biased random walk starting at $1$ does not visit zero at probability $(\lambda E\rho-1)/\lambda E\rho$.

For the classic SIS model with $\rho\equiv1$, Theorem \ref{theorem 2.1 main} shows that
\[
\lim_{d\rightarrow+\infty}P_{\lambda,d}\big(C_t^O\neq \emptyset,\forall~t\geq 0\big)=\frac{\lambda-1}\lambda
\]
for $\lambda\geq 1$. This result is first given in \cite{Xue2017} as far as we know.

Similar result with that in Theorem \ref{theorem 2.1 main} for the bond percolation model is obtained in \cite{Kesten1991}. In \cite{Kesten1991}, Kesten studies the high-dimensional Fortuin-Kasteleyn cluster model, containing the bond percolation model as a special case. It is shown in \cite{Kesten1991} that the probability that $O$ belongs to the infinite open cluster converges to the solution to the equation
\[
x=1-e^{-\lambda x}
\]
as the dimension $d$ grows to infinity for the bond percolation model on $\mathbb{Z}^d$ where an edge is open with probability $\frac{\lambda}{2d}$ with $\lambda>1$.

It is obviously that $P_{\lambda,d}\big(C_t^O\neq \emptyset,\forall~t\geq 0\big)$ and $P_{\lambda,d}\big(I_t^O\neq \emptyset,\forall~t\geq 0\big)$ are increasing with $\lambda$, hence it is reasonable to define
\[
\lambda_c(d)=\sup\big\{\lambda:~P_{\lambda,d}\big(C_t^O\neq \emptyset,\forall~t\geq 0\big)=0\big\}
\]
and
\[
\beta_c(d)=\sup\big\{\lambda:~P_{\lambda,d}\big(I_t^O\neq \emptyset,\forall~t\geq 0\big)=0\big\}.
\]
$\lambda_c(d)$ is called the critical value of the contact process, since it is the maximum of infection rates with which the infectious vertices die out with probability one. For similar reason, $\beta_c(d)$ is called the critical value of the SIR model.
The following conclusion about estimations of $\lambda_c(d)$ and $\beta_c(d)$ is an application of Theorem \ref{theorem 2.1 main}.
\begin{theorem}\label{theorem 2.2}
\[
\lim_{d\rightarrow+\infty}\beta_c(d)=\frac{1}{E\rho}\text{\quad while\quad}\limsup_{d\rightarrow+\infty}\lambda_c(d)\leq \frac{1}{E\rho}.
\]
\end{theorem}

When $\rho\equiv1$, Theorem \ref{theorem 2.2} shows that $\limsup_{d\rightarrow+\infty}\lambda_c(d)\leq 1$. A stronger conclusion that $\lim_{d\rightarrow+\infty}\lambda_c(d)=1$ for the classic contact process is proved by Holley and Liggett in \cite{Hol1981}. In \cite{Grif1983}, Griffeath gives another proof of this result and obtains a better upper bound of $\lambda_c(d)$. When $P(\rho=1)=p=1-P(\rho=0)$, Theorem \ref{theorem 2.2} shows that $\lim_{d\rightarrow+\infty}\beta_c(d)=\frac{1}{p}$ while $\limsup_{d\rightarrow+\infty}\lambda_c(d)\leq \frac{1}{p}$ for SIR and SIS models on clusters of bond percolation model. These two results are proved in \cite{Xue2017b} and \cite{Xue2016} respectively.

We believe that $\lim_{d\rightarrow+\infty}\lambda_c(d)=\frac{1}{E\rho}$ but have not found a proof yet. Since $P_{\lambda,d}\big(C_t^O\neq \emptyset,\forall~t\geq 0\big)$ is increasing with $\lambda$, a direct corollary of Theorem \ref{theorem 2.1 main} is that
\[
\lim_{d\rightarrow+\infty}P_{\lambda,d}\big(C_t^O\neq \emptyset,\forall~t\geq 0\big)=0
\]
for any $\lambda<\frac{1}{E\rho}$. If this conclusion can be strengthened to that
\[
P_{\lambda,d}\big(C_t^O\neq \emptyset,\forall~t\geq 0\big)=0
\]
for $\lambda<\frac{1}{E\rho}$ and sufficiently large $d$, then we can claim that $\liminf_{d\rightarrow+\infty}\lambda_c(d)\geq \frac{1}{E\rho}$ and hence
$\lim_{d\rightarrow+\infty}\lambda_c(d)=\frac{1}{E\rho}$. We will work on this problem as a further study.

According to the basic coupling of Markov process (see Section 3.1 of \cite{Lig1985}), it is easy to check that
\[
P_{\lambda,d}\big(C_t^O\neq \emptyset,\forall~t\geq 0\big)\geq P_{\lambda,d}\big(I_t^O\neq \emptyset,\forall~t\geq 0\big).
\]
Therefore, to prove Theorem \ref{theorem 2.1 main}, we only need to show that
\begin{equation}\label{equ 2.1}
\limsup_{d\rightarrow+\infty}P_{\lambda,d}\big(C_t^O\neq \emptyset,\forall~t\geq 0\big)\leq \frac{\lambda E\rho-1}{\lambda E\rho}
\end{equation}
and
\begin{equation}\label{equ 2.2}
\liminf_{d\rightarrow+\infty}P_{\lambda,d}\big(I_t^O\neq \emptyset,\forall~t\geq 0\big)\geq \frac{\lambda E\rho-1}{\lambda E\rho}
\end{equation}
for $\lambda>\frac{1}{E\rho}$.

The proof of Equation \eqref{equ 2.1} is given in Section \ref{section 3}. The core idea of the proof is as follows. For given large inter $M$ and small positive constant $\epsilon$, with high probability that every vertex $x$ in the set $\big\{u:~\|u\|\leq M\big\}$ satisfies that
\[
\sum_{y:y\sim x}\rho(x,y)\leq 2d (E\rho+\epsilon).
\]
Before the first moment when $C_t$ contains a vertex with $l_1$-norm larger than $M$, the embedded chain of $|C_t|$ is dominated from above by a biased random walk which increases by one with probability $\lambda(E\rho+\epsilon)/1+\lambda(E\rho+\epsilon)$ or decrease by one with probability $1/1+\lambda(E\rho+\epsilon)$. Such a biased random walk starting at $1$ hits zero at least once with probability $1/\lambda(E\rho+\epsilon)$.

The proof of Equation \eqref{equ 2.2} is given in Section \ref{section 4}. The core idea of the proof is as follows. We divide $\mathbb{Z}^d$ into two disjoint parts $\Gamma_1$ and $\Gamma_2$. We first show that there exist $d^{1/3}$ vertices in $\Gamma_1$ which are infected through paths on $\Gamma_1$ with probability about $(\lambda E\rho-1)/\lambda E\rho$. In this step we dominate the embedded chain of $|I_t|$ from below by another biased random walk. Then we show that these $d^{1/3}$ vertices infect at least $d^{1/4}$ vertices in $\Gamma_2$ by edges connecting $\Gamma_1$ and $\Gamma_2$ with high probability. At last, we show that with $d^{1/4}$ initial infectious vertices in $\Gamma_2$, the SIR model confined to $\Gamma_2$ survives with high probability. The approach in this step is inspired by the technique introduced in \cite{Xue2017b}. Since $\Gamma_1$ and $\Gamma_2$ are disjoint, the event concerned with in the third step is independent of the two events concerned with in the first and second steps and hence the survival probability is at least the product of the probability of the third event and the probability that both the first and the second events occur. To make the above explanation rigorous, we introduce the definition of so-called infectious path at the beginning of Section \ref{section 4} that a vertex $x$ has ever been infected when and only when there exists an infectious path from $O$ to $x$.

The proof of Theorem \ref{theorem 2.2} is given in Section \ref{section 5}, which is an application of Theorem \ref{theorem 2.1 main} and the definition of infectious path introduced in Section \ref{section 4}.

\section{Proof of Equation \eqref{equ 2.1}}\label{section 3}
In this section we give the proof of Equation \eqref{equ 2.1}. Throughout this section we assume that $\lambda\geq 1/E\rho$. First we introduce some notations and definitions. For $r>0$, we define
\[
B(d,r)=\big\{x\in \mathbb{Z}^d:~\|x\|\leq r\big\}
\]
as the set of vertices with $l_1$ norm at most $r$. For $M>0$ and $\epsilon>0$, we define
\[
A(d,M,\epsilon)=\big\{\omega\in \Omega_d:~\sum_{y:y\sim x}\rho(x,y,\omega)\leq 2d(E\rho+\epsilon),\forall~x\in B(d,M)\big\}
\]
as the set of random environments where every vertex $x$ with $l_1$ norm at most $M$ satisfies that
\[
\sum_{y:y\sim x}\rho(x,y)\leq 2d(E\rho+\epsilon).
\]
According to the classic theory about large deviation principle, there exists $J(\epsilon)>0$ that
\[
\mu_d\big(\omega:\frac{1}{2d}\sum_{y:y\sim x}\rho(x,y,\omega)\geq E\rho+\epsilon\big)\leq e^{-2dJ(\epsilon)}
\]
for $d\geq 1$ and given $x\in \mathbb{Z}^d$, since $\{\rho(x,y):~y\sim x\}$ are $2d$ independent copies of $\rho$. For each $x\in B(d,M)$, there is an path from $O$ to $x$ with length at most $M$. For a path on $\mathbb{Z}^d$, each step has $2d$ choices, therefore
\[
|B(d,M)|\leq \sum_{l=0}^M(2d)^l\leq (M+1)(2d)^M.
\]
As a result,
\begin{equation}\label{equ 3.1}
\mu_d\big(A(d,M,\epsilon)\big)\geq 1-(M+1)(2d)^Me^{-2dJ(\epsilon)}
\end{equation}
according to the Chebyshev's inequality.

We define $\{V_n\}_{n\geq 0}$ as the biased random walk on $\mathbb{Z}^1$ that
\[
P(V_{n+1}-V_n=1)=\frac{\lambda(E\rho+\epsilon)}{1+\lambda(E\rho+\epsilon)}=1-P(V_{n+1}-V_n=-1)
\]
for $n\geq 0$ and $V_0=1$. For $K\geq 0$, we define
\[
\tau_K=\inf\{n\geq 0:~V_n=K\}
\]
as the first moment when $K$ is visited. According to classic theory about biased random walk,
\begin{equation}\label{equ 3.2}
\lim_{K\rightarrow+\infty}\frac{\tau_K}{K}=\frac{1}{\frac{\lambda(E\rho+\epsilon)}{1+\lambda(E\rho+\epsilon)}-\frac{1}{1+\lambda(E\rho+\epsilon)}}
=\frac{1+\lambda(E\rho+\epsilon)}{\lambda(E\rho+\epsilon)-1}.
\end{equation}
Now we give the proof of Equation \eqref{equ 2.1}.

\proof[Proof of Equation \eqref{equ 2.1}]

For given $M>0$, $\epsilon>0$ and $\omega\in A(d,M,\epsilon)$, $C_t$ with edge weights $\{\rho(e,\omega)\}_{e\in E_d}$ decreases by one at rate $|C_t|$ or increases by one at rate at most
\[
\frac{\lambda}{2d}\sum_{x\in C_t}\sum_{y:y\sim x}\rho(x,y,\omega)\leq \frac{\lambda}{2d}2d(E\rho+\epsilon)|C_t|=\lambda(E\rho+\epsilon)|C_t|
\]
for $t<\inf\{s:~C_s\ni x\text{~for some~}x\text{~with~}\|x\|>M\}$. As a result, before the moment $\inf\{s:~C_s\ni x\text{~for some~}x\text{~with~}\|x\|>M\}$, the embedded chain of $|C_t|$ is dominated from above by $\{V_n\}_{n\geq 0}$. Since all the infections occur between nearest neighbors, the state of $\{C_t\}_{t\geq 0}$ must jump at least $M$ times to make $C_t$ contain a vertex with $l_1$ norm lager than $M$. According to the above analysis, for given $K>0$,
\[
\big\{C_t=\emptyset\text{~for some~}t>0\big\}\supseteq \{\tau_0<\tau_K,\tau_K<M\}
\]
in the sense of coupling, where $\tau_K$ is the first time $K$ is visited by $\{V_n\}_{n\geq 0}$ as we have defined. Therefore, for $\omega\in A(d,M,\epsilon)$,
\begin{align}\label{equ 3.3}
P_{\lambda,\omega}\big(C_t=\emptyset\text{~for some~}t>0\big)&\geq P(\tau_0<\tau_K)-P(\tau_K\geq M)\\
&=\frac{\frac{1}{\lambda(E\rho+\epsilon)}-(\frac{1}{\lambda(E\rho+\epsilon)})^K}{1-(\frac{1}{\lambda(E\rho+\epsilon)})^K}-P(\tau_K\geq M) \notag
\end{align}
according to classic theorem of biased random walk. Then, according to Equation \eqref{equ 3.1},
\begin{align}\label{equ 3.4}
&P_{\lambda,d}\big(C_t=\emptyset\text{~for some~}t>0\big)\notag\\
&\geq E_{\mu_d}\Big[P_{\lambda,\omega}\big(C_t=\emptyset\text{~for some~}t>0\big)1_{A(d,M,\epsilon)}\Big] \\
&\geq \Big[\frac{\frac{1}{\lambda(E\rho+\epsilon)}-(\frac{1}{\lambda(E\rho+\epsilon)})^K}{1-(\frac{1}{\lambda(E\rho+\epsilon)})^K}-P(\tau_K\geq M)\Big]\Big(1-(M+1)(2d)^Me^{-2dJ(\epsilon)}\Big).\notag
\end{align}
We choose $K=\lfloor d^{1/2}\rfloor$ and $M=2K\frac{1+\lambda(E\rho+\epsilon)}{\lambda(E\rho+\epsilon)-1}$, then $\lim_{d\rightarrow+\infty}P(\tau_K\geq M)=0$ by Equation \eqref{equ 3.2}, $\lim_{d\rightarrow+\infty}(M+1)(2d)^Me^{-2dJ(\epsilon)}=0$ and
\[
\lim_{d\rightarrow+\infty}\frac{\frac{1}{\lambda(E\rho+\epsilon)}-(\frac{1}{\lambda(E\rho+\epsilon)})^K}{1-(\frac{1}{\lambda(E\rho+\epsilon)})^K}
=\frac{1}{\lambda(E\rho+\epsilon)}.
\]
As a result, by Equation \eqref{equ 3.4},
\[
\liminf_{d\rightarrow+\infty}P_{\lambda,d}\big(C_t=\emptyset\text{~for some~}t>0\big)\geq \frac{1}{\lambda(E\rho+\epsilon)}.
\]
Since $\epsilon$ is arbitrary, we have
\begin{equation}\label{equ 3.5}
\liminf_{d\rightarrow+\infty}P_{\lambda,d}\big(C_t=\emptyset\text{~for some~}t>0\big)\geq \frac{1}{\lambda E\rho}.
\end{equation}
Equation \eqref{equ 2.1} follows from Equation \eqref{equ 3.5} directly.

\qed

\section{Proof of Equation \eqref{equ 2.2}}\label{section 4}
The aim of this section is to prove Equation \eqref{equ 2.2}. Throughout this section we assume that $\lambda>\frac{1}{E\rho}$, since the case where $\lambda=\frac{1}{E\rho}$ becomes trivial after the case where $\lambda>\frac{1}{E\rho}$ is proved. This section is divided into four parts. In Subsection \ref{subsection 4.0} we give the proof of Equation \eqref{equ 2.2} based on Lemmas \ref{lemma 4.3} and \ref{lemma 4.4}. The proof of Lemma \ref{lemma 4.4} is given in Subsection \ref{subsection 4.3} while the proof of Lemma \ref{lemma 4.3} is given in Subsection \ref{subsection 4.2}. The proof of Lemma \ref{lemma 4.3} utilizes Lemma \ref{lemma 4.2}. We give the proof of Lemma \ref{lemma 4.2} in Subsection \ref{subsection 4.1}.

\subsection{Proof of Equation \eqref{equ 2.2}}\label{subsection 4.0}
In this subsection we give the proof of Equation \eqref{equ 2.2}. First we introduce the definition of the infectious path. Let
\[
H_d=\big\{(x,y)\in \mathbb{Z}^d\times \mathbb{Z}^d:~x\sim y\big\}
\]
be the set of ordered pairs of neighbors on $\mathbb{Z}^d$, then we define
\[
X_3=[0,+\infty)^{\mathbb{Z}^d}\times[0,+\infty)^{H_d}.
\]
Therefore, an element in $X_3$ can be written as $(Y,U)$, where $Y: \mathbb{Z}^d\rightarrow [0,+\infty)$ and $U: H_d\rightarrow [0,+\infty)$. Let $\mathcal{F}_3$ be the smallest sigma-field that $\{Y(x)\}_{x\in \mathbb{Z}^d}$ and $\{U(x,y)\}_{(x,y)\in H_d}$ are measurable with respect to.

For any $\omega\in \Omega_d$, let $\nu_\omega$ be a probability measure on $(X_3,\mathcal{F}_3)$ that $Y(x)$ is an exponential time with rate one for each $x\in \mathbb{Z}^d$ and $U(y,z)$ is an exponential time with rate $\frac{\lambda}{2d}\rho(y,z,\omega)$ for each $(y,z)\in H_d$ while all these exponential times are independent under $\nu_\omega$.

For a self-avoiding path $\vec{l}=(l_0,l_1,\ldots,l_n)$ on $\mathbb{Z}^d$ with length $n$, we say $\vec{l}$ is an infectious path (with respect to $(Y,U)$) when and only when $U(l_i,l_{i+1})<Y(l_i)$ for $0\leq i\leq n-1$. We have the following important lemma.
\begin{lemma}\label{lemma 4.1}
Let $X_2$ be defined as in Section \ref{section 1} and $\mathcal{B}$ be the smallest sigma-field containing all the finite cylinder sets included in $X_2^{[0,+\infty)}$, then there exists a measurable mapping
\[
\{(\widehat{S}_t,\widehat{I}_t)\}_{t\geq 0}: (X_3,\mathcal{F}_3)\rightarrow \Big(X_2^{[0,+\infty)}, \mathcal{B}\Big)
\]
that $\{(\widehat{S}_t,\widehat{I}_t)\}_{t\geq 0}$ under the measure $\nu_\omega$ is a version of $\{(S_t^O,I_t^O)\}_{t\geq 0}$ under the probability measure $P_{\lambda,\omega}(\cdot)$ for each $\omega\in \Omega_d$ and
\begin{align*}
\Big\{(Y,U):& \text{~}x\in \widehat{I}_t\text{~for some~}t>0\Big\}=\Big\{(Y,U):\\
&\text{~there exists an infectious path with respect to $(Y,U)$ from $O$ to $x$}\Big\}
\end{align*}
for any $x\neq O$.
\end{lemma}

We omit the proof of Lemma \ref{lemma 4.1} here since it is a little tedious while this lemma can be explained intuitively and clearly. The intuitive explanation of Lemma \ref{lemma 4.1} is as follows. $Y(x)$ is the time $x$ waits for to become recovered after $x$ is infected, i.e., $x$ becomes recovered at moment $t+Y(x)$ if $x$ is infected at moment $t$. $U(x,y)$ is the time $x$ waits for to infect neighbor $y$ after $x$ is infected. The infection really occurs when $U(x,y)<Y(x)$ and $y$ is not infected by other vertices before the moment $t+U(x,y)$, where $t$ is the moment when $x$ is infected. As a result, for any $x\neq O$, if $x$ has ever been infected, then there exists a self-avoiding path $\vec{l}=(O,l_1,\ldots,l_{n-1},x)$ that $l_i$ has ever infected $l_{i+1}$ for $0\leq i\leq n-1$ and hence $U(l_i,l_{i+1})<Y(l_i)$ for $0\leq i\leq n-1$, i.e., $\vec{l}$ is an infectious path. On the other hand, if $\vec{l}=(O,l_1,\ldots,l_{n-1},x)$ is an infectious path, then we claim that $l_i$ has ever been infected for all $0\leq i\leq n$. This claim holds for $i=0$ trivially since $O\in I_0$. Assuming that $l_i$ has ever been infected for some $i<n$, then there are two possible cases. The first case is that $l_{i+1}\in I_s$ for some $s<\inf\{t:~l_i\in I_t\}+U(l_i,l_{i+1})$, then our claim holds for $i+1$ trivially. The second case is that $l_{i+1}\not\in I_s$ for any $s<\inf\{t:~l_i\in I_t\}+U(l_i,l_{i+1})$, then
\[
\inf\{t:l_{i+1}\in I_t\}=\inf\{t:~l_i\in I_t\}+U(l_i,l_{i+1})
\]
as we have introduced. As a result, our claim holds for $i+1$ and then holds for all $0\leq i\leq n$ according to the principle of mathematical induction. In conclusion,
\[
\big\{x\text{~has ever been infected}\big\}=\big\{\text{there is an infectious path from $O$ to $x$}\big\}.
\]

For simplicity, from now on we identify $\{(\widehat{S}_t,\widehat{I}_t)\}_{t\geq 0}$ given by Lemma \ref{lemma 4.1} with $\{(S_t^O,I_t^O)\}_{t\geq 0}$ and identify $\nu_\omega$ with $P_{\lambda,\omega}$. This identification is permitted by Lemma \ref{lemma 4.1}. As a result, $P_{\lambda,d}$ can be considered as a probability measure on $(X_3,\mathcal{F}_3)$ and
\[
P_{\lambda,d}(\cdot)=E_{\mu_d}\Big(P_{\lambda,\omega}(\cdot)\Big)=E_{\mu_d}\Big(\nu_\omega(\cdot)\Big).
\]
For later use, we introduce some definitions. We define
\[
\Gamma_1=\Big\{x=(x_1,\ldots,x_d)\in \mathbb{Z}^d:\sum_{i=d-\lfloor\frac{d}{\log d}\rfloor+1}^d|x_i|=0\Big\},
\]
\begin{align*}
\Gamma_2=\Big\{x=(x_1,\ldots,x_d)\in \mathbb{Z}^d:~&x_i\geq 0 \text{~for all~}d-\lfloor\frac{d}{\log d}\rfloor+1\leq i\leq d \\ &\text{~and~}\sum_{i=d-\lfloor\frac{d}{\log d}\rfloor+1}^dx_i>0\Big\}
\end{align*}
and
\[
\Gamma_3=\Big\{x=(x_1,\ldots,x_d)\in \Gamma_2:\sum_{i=d-\lfloor\frac{d}{\log d}\rfloor+1}^d x_i=1\Big\},
\]
where $\lfloor u\rfloor=n$ when $n$ is an integer and $n\leq u<n+1$.

We say an infectious path is on a subgraph $A$ of $\mathbb{Z}^d$ when all the vertices on this path belong to $A$.  We define
\[
D_1=\Big\{x\in \Gamma_1:\text{~there is an infectious path on $\Gamma_1$ from $O$ to $x$}\Big\}.
\]
Note that $D_1$ is a mapping from $X_3$ to the power set of $\Gamma_1$. The following lemma is important for us to prove Equation \eqref{equ 2.2}.
\begin{lemma}\label{lemma 4.2}
\[
\liminf_{d\rightarrow+\infty}P_{\lambda,d}\Big(|D_1|\geq \lfloor d^{1/3}\rfloor\Big)\geq \frac{\lambda E\rho-1}{\lambda E\rho}.
\]
\end{lemma}

The proof of Lemma \ref{lemma 4.2} will be given in Subsection \ref{subsection 4.1}.

For any $x\in \Gamma_1$ and any $B\subseteq \Gamma_1$, we define
\[
D_2(x)=\Big\{y\in \Gamma_3:~y\sim x \text{~and~}U(x,y)<Y(x)\Big\}
\]
and $D_2(B)=\bigcup_{w\in B}D_2(w)$. The following lemma about $D_2(D_1)$ is important for us to prove Equation \eqref{equ 2.2}.
\begin{lemma}\label{lemma 4.3}
\[
\liminf_{d\rightarrow+\infty}P_{\lambda,d}\big(|D_2(D_1)|>d^{1/4}\big)\geq \frac{\lambda E\rho-1}{\lambda E\rho}.
\]
\end{lemma}
The proof of Lemma \ref{lemma 4.3} is given in Subsection \ref{subsection 4.2}, where Lemma \ref{lemma 4.2} will be utilized.

For integer $n\geq 1$ and $B\subseteq \Gamma_3$, we define
\begin{align*}
D_3(n,B)=\Big\{x\in \Gamma_2:&~\|x\|\geq n\text{~and there exists an infectious path}\\
&\text{on $\Gamma_2$ from some vertex in $B$ to~}\|x\|\Big\}.
\end{align*}
Note that $D_3(n,B)$ is a mapping from $X_3$ to the power set of $\Gamma_2$ for given $n$ and $B$. The following lemma is crucial for us to prove Equation \eqref{equ 2.2}, where we use $\{A_n \text{~i.o.~}\}$ to denote $\bigcap_{n\geq 1}\bigcup_{k\geq n} A_k$ for a series of events $\{A_n\}_{n\geq 1}$.

\begin{lemma}\label{lemma 4.4}
For each $d\geq 1$, let
\[
\Delta(d)=\inf\Big\{P_{\lambda,d}\big(D_3(n,B)\neq \emptyset \text{~i.o.}\big):~B\subseteq \Gamma_3\text{~and~}|B|=\lfloor d^{1/4}\rfloor\Big\},
\]
then
\[
\lim_{d\rightarrow+\infty}\Delta(d)=1.
\]
\end{lemma}
The proof of Lemma \ref{lemma 4.4} is given in Subsection \ref{subsection 4.3}. The strategy of the proof is inspired by the approach introduced in \cite{Xue2017b}.

At the end of this subsection we show how to utilize Lemmas \ref{lemma 4.3} and \ref{lemma 4.4} to prove Equation \eqref{equ 2.2}.

\proof[Proof of Equation \eqref{equ 2.2}]

According to the definitions of $D_1, D_2(B)$ and $D_3(n,B)$, for each $n\geq 1$ and any $x\in D_3(n,D_2(D_1))$, there exist $y\in \Gamma_1$ and $z\in \Gamma_3$ that the following three conditions holds.

(1) There is an infectious path on $\Gamma_2$ from $z$ to $x$.

(2) $y\sim z$ and $U(y,z)<Y(y)$.

(3) There is an infectious path on $\Gamma_1$ from $O$ to $y$.

\begin{figure}[H]
  \centering
  \includegraphics[height=0.4\textwidth=0.4]{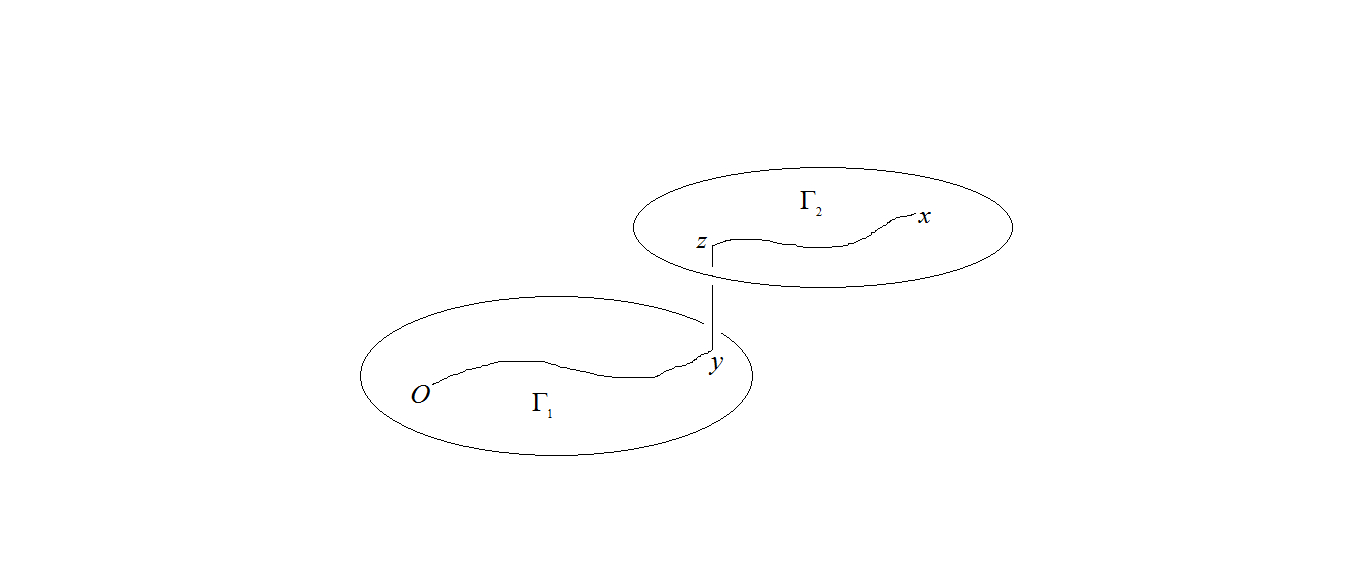}
  \caption{Infectious path}
  \label{graph 4.1}
\end{figure}

As a result, there is an infectious path from $O$ to $x$, as it is shown in Figure \ref{graph 4.1}.

Therefore, by Lemma \ref{lemma 4.1},
\begin{equation*}
D_3\big(n,D_2(D_1)\big)\subseteq \bigcup_{t\geq 0}I_t
\end{equation*}
and hence
\begin{equation*}
\Big\{D_3\big(n,D_2(D_1)\big)\neq \emptyset \text{\quad i.o.~}\Big\}\subseteq \Big\{|\bigcup_{t\geq 0}I_t|=+\infty\Big\}.
\end{equation*}
If there are infinite many vertices have ever been infected, then they can not all become recovered before a uniform moment $T<+\infty$, since each infected vertex waits for an independent copy of the exponential time with rate $1$ to become recovered. As a result,
\[
\Big\{|\bigcup_{t\geq 0}I_t|=+\infty\Big\}\subseteq \big\{I_t\neq \emptyset,\forall~t\geq 0\big\}
\]
and hence
\begin{equation}\label{equ 40.1}
P_{\lambda,d}\big(I_t\neq \emptyset,\forall~t\geq 0\big)\geq P_{\lambda,d}\Big(D_3\big(n,D_2(D_1)\big)\neq \emptyset \text{\quad i.o.~}\Big).
\end{equation}
Conditioned on $D_2(D_1)=B$ for some $B\subseteq \Gamma_3$, $D_3(n,B)$ depends on $\{Y(x):~x\in \Gamma_2\}\bigcup\{U(y,z):~y,z\in \Gamma_2\}$, which is independent of the random set $D_2(D_1)$, since $D_2(D_1)$ depends on $\{Y(x):~x\in \Gamma_1\}\bigcup\{U(y,z):~y\in \Gamma_1,z\in \Gamma_1\bigcup\Gamma_3\}$. Therefore,
\begin{align}\label{equ 40.2}
&P_{\lambda,d}\Big(D_3\big(n,D_2(D_1)\big)\neq \emptyset \text{\quad i.o.~}\Big|D_2(D_1)\Big)\notag\\
&=P_{\lambda,d}\Big(D_3\big(n,B\big)\neq \emptyset \text{\quad i.o.~}\Big)\Bigg|_{B=D_2(D_1)}\\
&\geq \inf\Big\{P_{\lambda,d}\Big(D_3\big(n,A\big)\neq \emptyset \text{\quad i.o.~}\Big):~A\subseteq \Gamma_3\text{~and~}|A|=k\Big\}\Bigg|_{k=|D_2(D_1)|}.\notag
\end{align}
It is obviously that
\[
\inf\Big\{P_{\lambda,d}\Big(D_3\big(n,A\big)\neq \emptyset \text{\quad i.o.~}\Big):~A\subseteq \Gamma_3\text{~and~}|A|=k\Big\}
\]
is increasing with $k$. As a result, by Equation \eqref{equ 40.2},
\begin{align}\label{equ 40.3}
&P_{\lambda,d}\Big(D_3\big(n,D_2(D_1)\big)\neq \emptyset \text{\quad i.o.~}\Big)\notag\\
&=P_{\lambda,d}\Big(D_3\big(n,D_2(D_1)\big)\neq \emptyset \text{\quad i.o.~}\Big||D_2(D_1)|>\lfloor d^{1/4}\rfloor\Big)\\
&\text{\quad}\times P_{\lambda,d}\big(|D_2(D_1)|>\lfloor d^{1/4}\rfloor\big)\notag\\
&\geq E_{\lambda,d}\Bigg(\inf\Big\{P_{\lambda,d}\Big(D_3\big(n,A\big)\neq \emptyset \text{\quad i.o.~}\Big):~A\subseteq \Gamma_3\text{~and~}|A|=k\Big\}\Bigg|_{k=|D_2(D_1)|}\notag\\
&\text{\quad\quad\quad\quad}\Bigg||D_2(D_1)|>\lfloor d^{1/4}\rfloor\Bigg)\times P_{\lambda,d}\big(|D_2(D_1)|>\lfloor d^{1/4}\rfloor\big)\notag\\
&\geq \inf\Big\{P_{\lambda,d}\Big(D_3\big(n,A\big)\neq \emptyset \text{\quad i.o.~}\Big):~A\subseteq \Gamma_3\text{~and~}|A|=\lfloor d^{1/4}\rfloor\Big\}\notag\\
&\text{\quad}\times P_{\lambda,d}\big(|D_2(D_1)|>\lfloor d^{1/4}\rfloor\big)\notag\\
&=\Delta(d)P_{\lambda,d}\big(|D_2(D_1)|>\lfloor d^{1/4}\rfloor\big).\notag
\end{align}
By Equation \eqref{equ 40.3}, Lemmas \ref{lemma 4.3} and \ref{lemma 4.4},
\begin{equation}\label{equ 40.4}
\liminf_{d\rightarrow+\infty}P_{\lambda,d}\Big(D_3\big(n,D_2(D_1)\big)\neq \emptyset \text{\quad i.o.~}\Big)\geq \frac{\lambda E\rho-1}{\lambda E\rho}.
\end{equation}
Equation \eqref{equ 2.2} follows from Equations \eqref{equ 40.1} and \eqref{equ 40.4} directly.

\qed

\subsection{Proof of Lemma \ref{lemma 4.2}}\label{subsection 4.1}
In this subsection we give the proof of Lemma \ref{lemma 4.2}, which is similar with that of Equation \eqref{equ 2.1}. First we introduce some notations. For given $M>0$ and $\epsilon>0$, we define
\[
F(d,M,\epsilon)=\Big\{\omega\in \Omega:~\sum_{y:y\sim x}\rho(x,y,\omega)>2d(E\rho-\epsilon)\text{~for all~}x\in B(d,M)\Big\}.
\]
According to a similar analysis with that of Equation \eqref{equ 3.1}, there exist $J_2(\epsilon)>0$ such that
\begin{equation}\label{equ 41.1}
\mu_d\Big(F(d,M,\epsilon)\Big)\geq 1-(M+1)(2d)^Me^{-2dJ_2(\epsilon)}.
\end{equation}
We choose $\epsilon$ sufficiently small such that $\lambda(E\rho-2\epsilon)>1$. Then we assume that we deal with $d$ sufficiently large such that
\begin{equation}\label{equ 41.2}
E\rho-\epsilon-\frac{\big(1+M_2(d)\big)\Theta}{2d}-\frac{\Theta}{2d}\lfloor\frac{d}{\log d}\rfloor\geq E\rho-2\epsilon
\end{equation}
where $\Theta$ is defined as in Section \ref{section 1} while $M_2(d)=2\lfloor d^{1/3} \rfloor\frac{\lambda(E\rho-2\epsilon)+1}{\lambda(E\rho-2\epsilon)-1}$.
We define $\{W_n\}_{n\geq 0}$ as biased random walk on $\mathbb{Z}^1$ that $W_0=1$ and
\[
P(W_{n+1}-W_n=1)=\frac{\lambda(E\rho-2\epsilon)}{\lambda(E\rho-2\epsilon)+1}=1-P(W_{n+1}-W_n=-1).
\]
For each integer $K\geq 0$, we define
\[
\phi_K=\inf\{n\geq 0:~W_n=K\}
\]
as the first moment when $K$ is visited.

Now we give the proof of Lemma \ref{lemma 4.2}.
\proof[Proof of Lemma \ref{lemma 4.2}]

We denote by $\{(\widetilde{S}_t,\widetilde{I}_t):t\geq 0\}$ the SIR model with random edge weights confined to the graph $\Gamma_1$ with $(\widetilde{S}_0,\widetilde{I}_0)=(\Gamma_1\setminus\{O\},\{O\})$. Let $\widetilde{I}=\cup_{t\geq 0}\widetilde{I}_t$, then according to a similar analysis with that leads to Lemma \ref{lemma 4.1},
\begin{equation}\label{equ 41.3}
D_1=\widetilde{I}.
\end{equation}
Let $\psi=\inf\{t:|\widetilde{I}_t|=\lfloor d^{1/3}\rfloor\}$, then by Equation \eqref{equ 41.3},
\begin{equation}\label{equ 41.4}
P_{\lambda,d}\big(|D_1|\geq \lfloor d^{1/3}\rfloor\big)\geq P_{\lambda,d}\big(\psi<+\infty\big).
\end{equation}
When the number of jumps of the state of $\widetilde{I}_t$ is no more than $M_2(d)$, there are at most $1+M_2(d)$ vertices that have ever been infected and all the infectious vertices belong to $B(d,M_2(d))$. As a result, for $\omega\in F(d,M_2(d),\epsilon)$ and $\{(\widetilde{S}_t,\widetilde{I}_t):t\geq 0\}$ with random edge weights with respect to $\omega$, $|\widetilde{I}_t|$ decreases by one at rate $|\widetilde{I}_t|$ while increases by one at rate
\begin{align*}
\frac{\lambda}{2d}\sum_{x\in \widetilde{I}_t}\sum_{y:y\sim x}\rho(x,y)&-\frac{\lambda}{2d}\sum_{x\in \widetilde{I}_t}\sum_{y:y\sim x}\rho(x,y)1_{\{y\not\in \Gamma_1\}}\\
&-\frac{\lambda}{2d}\sum_{x\in \widetilde{I}_t}\sum_{y:y\sim x}\rho(x,y)1_{\{y\in \widetilde{I}_s\text{~for some~}s\geq t\}}\\
\geq \lambda|\widetilde{I}_t|(E\rho-\epsilon)&-\frac{\lambda \Theta}{2d}|\widetilde{I}_t|\lfloor\frac{2d}{\log d}\rfloor-\frac{\lambda}{2d}|\widetilde{I}_t|
(M_2(d)+1)\Theta\\
=\lambda |\widetilde{I}_t|\Big(E\rho-\epsilon&-\frac{\big(1+M_2(d)\big)\Theta}{2d}-\frac{\Theta}{2d}\lfloor\frac{d}{\log d}\rfloor\Big)\\
\geq \lambda |\widetilde{I}_t|(E\rho-2\epsilon&) \text{\quad}(\text{This step utilizes Equation \eqref{equ 41.2}.})
\end{align*}
before the moment when the state of $\{|\widetilde{I}_t|\}_{t\geq 0}$ jumps for the $M_2(d)$th time. As a result, the embedded chain of $\{|\widetilde{I}_t|\}_{t\geq 0}$ is dominated from below by $\{W_n\}_{n\geq 0}$ for $0\leq n<M_2(d)$. Therefore, for $\omega\in F(d,M_2(d),\epsilon)$ and $\{(\widetilde{S}_t,\widetilde{I}_t)\}_{t\geq 0}$ with random edge weights with respect to $\omega$,
\[
\Big\{\psi<+\infty\Big\}\supseteq \Big\{\phi_{\lfloor d^{1/3}\rfloor}<\phi_0,\phi_{\lfloor d^{1/3}\rfloor}<M_2(d)\Big\}
\]
in the sense of coupling. Therefore, for $\omega\in F(d,M_2(d),\epsilon)$,
\begin{align}\label{equ 41.5}
P_{\lambda,\omega}\big(\psi<\infty\big)&\geq P\big(\phi_{\lfloor d^{1/3}\rfloor}<\phi_0,\phi_{\lfloor d^{1/3}\rfloor}<M_2(d)\big)\notag\\
&\geq  P\big(\phi_{\lfloor d^{1/3}\rfloor}<\phi_0\big)-P\big(\phi_{\lfloor d^{1/3}\rfloor}\geq M_2(d)\big)\\
&= \frac{1-\frac{1}{\lambda(E\rho-2\epsilon)}}{1-(\frac{1}{\lambda(E\rho-2\epsilon)})^{\lfloor d^{1/3}\rfloor}}-P\big(\phi_{\lfloor d^{1/3}\rfloor}\geq M_2(d)\big)\notag
\end{align}
according to the classic theory of biased random walk. By Equations \eqref{equ 41.1} and \eqref{equ 41.5},
\begin{align}\label{equ 41.6}
P_{\lambda,d}\big(\psi<\infty\big)\geq &E_{\mu_d}\Big[P_{\lambda,\omega}\big(\psi<\infty\big)1_{F(d,M_2(d),\epsilon)}\Big]\notag\\
\geq &\Big(\frac{1-\frac{1}{\lambda(E\rho-2\epsilon)}}{1-(\frac{1}{\lambda(E\rho-2\epsilon)})^{\lfloor d^{1/3}\rfloor}}-P\big(\phi_{\lfloor d^{1/3}\rfloor}\geq M_2(d)\big)\Big)\\
&\times\Big(1-(M_2(d)+1)(2d)^{M_2(d)}e^{-2dJ_2(\epsilon)}\Big). \notag
\end{align}
According to a similar analysis with that of Equation \eqref{equ 3.2},
\[
\lim_{d\rightarrow+\infty}\frac{\phi_{\lfloor d^{1/3}\rfloor}}{M_2(d)}=\frac{1}{2}
\]
and hence
\begin{equation}\label{equ 41.7}
\lim_{d\rightarrow+\infty}P\big(\phi_{\lfloor d^{1/3}\rfloor}\geq M_2(d)\big)=0.
\end{equation}
Then by Equations \eqref{equ 41.6} and \eqref{equ 41.7},
\begin{equation*}
\liminf_{d\rightarrow+\infty}P_{\lambda,d}\big(\psi<\infty\big)\geq 1-\frac{1}{\lambda(E\rho-2\epsilon)}
\end{equation*}
since
\[
\lim_{d\rightarrow+\infty}1-(M_2(d)+1)(2d)^{M_2(d)}e^{-2dJ_2(\epsilon)}=1.
\]
Since $\epsilon$ is arbitrary, we have
\begin{equation}\label{equ 41.8}
\liminf_{d\rightarrow+\infty}P_{\lambda,d}\big(\psi<\infty\big)\geq 1-\frac{1}{\lambda E\rho}=\frac{\lambda E\rho-1}{\lambda E\rho}.
\end{equation}
Lemma \ref{lemma 4.2} follows from Equations \eqref{equ 41.4} and \eqref{equ 41.8} directly.

\qed

\subsection{Proof of Lemma \ref{lemma 4.3}}\label{subsection 4.2}
In this subsection we give the proof of Lemma \ref{lemma 4.3}. First we introduce some notations and definitions. Let $\{\Psi(x)\}_{x\in \mathbb{Z}^d}$ be i.i.d. exponential times with rate $\lambda\Theta$ and independent with $\{Y(x)\}_{x\in \mathbb{Z}^d}$ and $\{U(x,y)\}_{(x,y)\in H_d}$ under the measure $P_{\lambda,\omega}$ for any $\omega\in \Omega$, where $\Theta$ is defined as in Section \ref{section 1}. Note that to make the above definition rigorous we can expand $X_3$ to $\widetilde{X}_3=X_3\times [0,+\infty)^{\mathbb{Z}^d}$ and identify $P_{\lambda,\omega}$ with the measure $\nu_\omega\times \pi$, where $\pi$ is the probability measure of i.i.d exponential times with rate $\lambda\Theta$. This is classic approach in measure theory so we omit the details.

For any $A\subseteq \Gamma_1$, we denote by $q(A)$ the random event that $Y(x)<\Psi(x)$ for any $x\in A$. For any $s>0$ and $A\subseteq \Gamma_1$, it is easy to check that \[E_{\lambda,d}\Big(e^{-s|D_2(A)|}\Big|q(A)\Big)\] depends only on $s$ and the cardinality of $A$. Hence we can reasonably define
\[
h(d,s,K)=E_{\lambda,d}\Big(e^{-s|D_2(A)|}\Big|q(A)\Big)
\]
for $s>0$ and $A\subseteq \Gamma_1$ with $|A|=K$. The following lemma is crucial for us to prove Lemma \ref{lemma 4.3}.

\begin{lemma}\label{lemma 42.1}
For any $s>0$,
\[
\lim_{d\rightarrow+\infty}h(d,-\frac{\log d}{d^{1/3}}s,\lfloor d^{1/3}\rfloor)=\exp\{-\frac{\lambda sE\rho}{2(\lambda \Theta+1)}\}.
\]
\end{lemma}

The intuitive explanation of Lemma \ref{lemma 42.1} is as follows. By direct calculation, it is easy to check that
\[
E_{\lambda,d}\Big(D_2(A)\Big|q(A)\Big)\approx \frac{\lambda d^{1/3}E\rho}{2\log d(\lambda \Theta+1)}
\]
for $A\subseteq \Gamma_1$ with $|A|=\lfloor d^{1/3}\rfloor$ and large $d$. Then it is natural to check wether
the distribution of $\frac{\log d}{d^{1/3}}|D_2(A)|$ conditioned on $q(A)$ converges weakly to the Dirac measure on $\frac{\lambda E\rho}{2(\lambda \Theta+1)}$.
One approach to do so is the Laplace transform, i.e., the calculation of $h(d,-\frac{\log d}{d^{1/3}}s,\lfloor d^{1/3}\rfloor)$.

We give the proof of Lemma \ref{lemma 42.1} at the end of this subsection. Now we show how to utilize Lemmas \ref{lemma 4.2} and \ref{lemma 42.1} to prove Lemma \ref{lemma 4.3}.

\proof[Proof of Lemma \ref{lemma 4.3}]

By Lemma \ref{lemma 4.2},
\begin{align}\label{equ 42.9}
&\liminf_{d\rightarrow+\infty}P_{\lambda,d}\big(|D_2(D_1)|>d^{1/4}\big)\notag\\
&\geq \liminf_{d\rightarrow+\infty}P_{\lambda,d}\Big(|D_2(D_1)|>d^{1/4}\Big||D_1|\geq \lfloor d^{1/3}\rfloor\Big)P_{\lambda,d}\Big(|D_1|\geq \lfloor d^{1/3}\rfloor\Big) \\
&\geq \frac{\lambda E\rho-1}{\lambda E\rho}\liminf_{d\rightarrow+\infty}P_{\lambda,d}\Big(|D_2(D_1)|>d^{1/4}\Big||D_1|\geq \lfloor d^{1/3}\rfloor\Big). \notag
\end{align}
We claim that
\begin{equation}\label{equ 42.10}
P_{\lambda,d}\Big(|D_2(D_1)|>d^{1/4}\Big||D_1|\geq \lfloor d^{1/3}\rfloor\Big)\geq P_{\lambda,d}\Big(|D_2(A)|>d^{1/4}\Big|q(A)\Big),
\end{equation}
where $|A|=\lfloor d^{1/3}\rfloor$. Note that $P_{\lambda,d}\Big(|D_2(A)|>d^{1/4}\Big|q(A)\Big)$ depends only on $|A|$, not on the choice of $A$. The explanation of Equation \eqref{equ 42.10} is as follows. Conditioned on $|D_1|\geq \lfloor d^{1/3}\rfloor$, there exists a subset $B$ of $D_1$ that $|B|=\lfloor d^{1/3}\rfloor$. If $|D_2(B)|>d^{1/4}$, then $|D_2(D_1)|>d^{1/4}$. The event $|D_1|\geq \lfloor d^{1/3}\rfloor$ relies on the values of $\{Y(x)\}_{x\in \Gamma_1}$ and $\{U(x,y)\}_{x\sim y,x,y\in \Gamma_1}$. So the event $|D_1|\geq \lfloor d^{1/3}\rfloor$ is correlated with the event $|D_2(B)|>d^{1/4}$ and we do not ensure (though we guess) that they are positive correlated . However, the worst condition with respect to $Y(\cdot)$ and $U(\cdot,\cdot)$ on $\Gamma_1$ for the probability that $|D_2(B)|>d^{1/4}$ occurs is that
\[
Y(x)< \inf\{U(x,y):~y\sim x, y\in \Gamma_1\}
\]
for any $x\in B$. Hence the probability that $|D_2(B)|>d^{1/4}$ occurs decreases if we replace the condition $|D_1|\geq \lfloor d^{1/3}\rfloor$ by that $Y(x)< \inf\{U(x,y):~y\sim x, y\in \Gamma_1\}$ for each $x\in B$. $\inf\{U(x,y):~y\sim x,y\in \Gamma_1\}$ is an exponential time with rate
\[
\sum_{y:y\sim x,y\in \Gamma_1}\frac{\lambda \rho(x,y)}{2d}\leq \frac{\lambda \Theta}{2d}(2d-2\lfloor\frac{d}{\log d}\rfloor)\leq \lambda \Theta.
\]
$\lambda\Theta$ is the rate of the exponential time $\Psi(x)$. As a result, the probability that $|D_2(B)|>d^{1/4}$ occurs will further decrease if we replace the condition $Y(x)< \inf\{U(x,y):~y\sim x, y\in \Gamma_1\}$ by $Y(x)<\Psi(x)$ for every $x\in B$, which leads to Equation \eqref{equ 42.10}.

For $A\subseteq \Gamma_1$ with $|A|=\lfloor d^{1/3}\rfloor$ and any $s>0$, by Chebyshev's inequality,
\begin{align}\label{equ 42.11}
P_{\lambda,d}\Big(|D_2(A)|<d^{1/4}\Big|q(A)\Big)&=P_{\lambda,d}\Big(e^{-s\frac{\log d}{d^{1/3}}|D_2(A)|}>e^{-s\frac{d^{1/4}}{d^{1/3}}\log d}\Big|q(A)\Big)\notag\\
&\leq \exp\big\{s\frac{d^{1/4}}{d^{1/3}}\log d\big\}h(d,-\frac{\log d}{d^{1/3}}s,\lfloor d^{1/3}\rfloor).
\end{align}
Then by Lemma \ref{lemma 42.1},
\begin{equation}\label{equ 42.12}
\limsup_{d\rightarrow+\infty}P_{\lambda,d}\Big(|D_2(A)|<d^{1/4}\Big|q(A)\Big) \leq \exp\{-\frac{\lambda sE\rho}{2(\lambda \Theta+1)}\}.
\end{equation}
Since $s$ is arbitrary, let $s\rightarrow+\infty$, we have
\[
\limsup_{d\rightarrow+\infty}P_{\lambda,d}\Big(|D_2(A)|<d^{1/4}\Big|q(A)\Big)=0
\]
and hence
\begin{equation}\label{equ 42.13}
\lim_{d\rightarrow+\infty}P_{\lambda,d}\Big(|D_2(A)|>d^{1/4}\Big|q(A)\Big)=1.
\end{equation}
Lemma \ref{lemma 4.3} follows from Equations \eqref{equ 42.9}, \eqref{equ 42.10} and \eqref{equ 42.13} directly.

\qed

At the end of this subsection we give the proof of Lemma \ref{lemma 42.1}.

\proof[Proof of Lemma \ref{lemma 42.1}]

According to assumptions of our model, it is easy to check that
\begin{equation}\label{equ 42.1}
h(d,-\frac{\log d}{d^{1/3}}s,\lfloor d^{1/3}\rfloor)=\Bigg(E_{\lambda,d}\Big(e^{-\frac{\log d}{d^{1/3}}s|D_2(O)|}\Big|Y(O)<\Psi(O)\Big)\Bigg)^{\lfloor d^{1/3}\rfloor}.
\end{equation}
and
\begin{align}\label{equ 42.2}
&E_{\lambda,d}\Big(e^{-\frac{\log d}{d^{1/3}}s|D_2(O)|}\Big|Y(O)<\Psi(O)\Big)\\
&=E_{\lambda,d}\Bigg(E_{\lambda,d}\Big(e^{-\frac{\log d}{d^{1/3}}s|D_2(O)|}\Big|Y(O)\Big)\Bigg|Y(O)<\Psi(O)\Bigg).\notag
\end{align}
Note that in Equation \eqref{equ 42.1} we utilize the fact that $|D_2(B)|=\sum_{x\in B}|D_2(x)|$ since $D_2(x)\cap D_2(y)=\emptyset$ for $x\neq y$.

According to assumptions of the model, it is easy to check that
\begin{equation}\label{equ 42.3}
E_{\lambda,d}\Big(e^{-\frac{\log d}{d^{1/3}}s|D_2(O)|}\Big|Y(O)\Big)=\Bigg[E_{\lambda,d}\Big(e^{-\frac{\log d}{d^{1/3}}s\Lambda_1}\Big|Y(O)\Big)\Bigg]^{\lfloor\frac{d}{\log d}\rfloor},
\end{equation}
where
\[
\Lambda_1=
\begin{cases}
1 & \text{~if~}U(O,e_d)<Y(O),\\
0 & \text{~if~}U(O,e_d)\geq Y(O)
\end{cases}
\]
and $e_d=(0,\ldots,0,1)$ as we have defined in Section \ref{section 1}.

By direct calculation,
\begin{align}\label{equ 42.4}
E_{\lambda,d}\Big(e^{-\frac{\log d}{d^{1/3}}s\Lambda_1}\Big|Y(O)\Big)
&=E\Big(e^{-\frac{\log d}{d^{1/3}}s}(1-e^{-\frac{\lambda\rho t}{2d}})+e^{-\frac{\lambda\rho t}{2d}}\Big)\Big|_{t=Y(O)} \\
&=\Big(1-\big(1-e^{-\frac{\log d}{d^{1/3}}s}\big)E\big(1-e^{-\frac{\lambda\rho t}{2d}}\big)\Big)\Big|_{t=Y(O)}, \notag
\end{align}
where $\rho$ is as defined in Section \ref{section 1} while $E$ is the expectation operator with respect to $\rho$. By Lagrange Mean Value Theorem and the fact that $e^a=1+a+o(a)$,
it is not difficult to check that
\begin{equation}\label{equ 42.5}
\lim_{d\rightarrow+\infty}\lfloor d^{1/3}\rfloor\Bigg\{\Big(1-\big(1-e^{-\frac{\log d}{d^{1/3}}s}\big)E\big(1-e^{-\frac{\lambda\rho t}{2d}}\big)\Big)^{\lfloor\frac{d}{\log d}\rfloor}-1\Bigg\}=-\frac{s\lambda tE\rho }{2}.
\end{equation}
By Equations \eqref{equ 42.3}, \eqref{equ 42.4} and \eqref{equ 42.5},
\begin{equation}\label{equ 42.6}
\lim_{d\rightarrow+\infty}\lfloor d^{1/3}\rfloor\Bigg\{E_{\lambda,d}\Big(e^{-\frac{\log d}{d^{1/3}}s|D_2(O)|}\Big|Y(O)\Big)-1\Bigg\}=-\frac{s\lambda Y(O)E\rho}{2}.
\end{equation}
By Equation \eqref{equ 42.6} and Dominated Convergence Theorem,
\begin{align}\label{equ 42.7}
&\lim_{d\rightarrow+\infty}\lfloor d^{1/3}\rfloor\Bigg\{E_{\lambda,d}\Big(e^{-\frac{\log d}{d^{1/3}}s|D_2(O)|}\Big|Y(O)<\Psi(O)\Big)-1\Bigg\}\\
&=-\frac{s\lambda E\rho}{2}E_{\lambda,d}\Big(Y(O)\Big|Y(O)<\Psi(O)\Big)=-\frac{\lambda sE\rho}{2(\lambda\Theta+1)}.\notag
\end{align}
According to classic conclusion about calculus, if $a_d\rightarrow 0$, $c_d\rightarrow +\infty$ while $a_dc_d\rightarrow b$ as $d\rightarrow+\infty$, then
\[
\lim_{d\rightarrow+\infty}(1+a_d)^{c_d}=e^b.
\]
Therefore, by Equation \eqref{equ 42.7},
\begin{equation}\label{equ 42.8}
\lim_{d\rightarrow+\infty} \Bigg(E_{\lambda,d}\Big(e^{-\frac{\log d}{d^{1/3}}s|D_2(O)|}\Big|Y(O)<\Psi(O)\Big)\Bigg)^{\lfloor d^{1/3}\rfloor}=\exp\{-\frac{\lambda sE\rho}{2(\lambda\Theta+1)}\}.
\end{equation}
Lemma \ref{lemma 42.1} follows from Equations \eqref{equ 42.1} and \eqref{equ 42.8} directly.

\qed

\subsection{Proof of Lemma \ref{lemma 4.4}}\label{subsection 4.3}
In this subsection we give the proof of Lemma \ref{lemma 4.4}. First we introduce some definitions and notations. For each $n\geq 1$, we use $\Xi_n$ to denote the set of self-avoiding paths on $\mathbb{Z}^d$ with length $n$.  For each $n\geq 1$, $x\in \Gamma_3$ and $B\subseteq \Gamma_3$, we define
\begin{align*}
L_n(x)=\Big\{\vec{x}=&(x_0,\ldots,x_n)\in \Xi_{n}:~x_0=x, \\
&x_{i}-x_{i-1}\in \big\{\pm e_j:~1\leq j\leq d-\lfloor\frac{d}{\log d}\rfloor\big\}
\text{~for each $i$ that~} \lfloor \log d\rfloor \nmid i \\
\text{~while~}&x_{i}-x_{i-1}\in \big\{e_j:~d-\lfloor\frac{d}{\log d}\rfloor+1\leq j\leq d\big\}\text{~for each $i$ that~}\lfloor\log d\rfloor\mid i
\Big\}
\end{align*}
and $L_n(B)=\bigcup_{w\in B}L_n(w)$, where we use $a\mid b$ to denote that $b$ is divisible by $a$ and $\{e_j:~1\leq j\leq d\}$ are defined as in Section \ref{section 1}.

For any $\vec{x},\vec{y}\in L_n(B)$, we define
\[
\sigma(\vec{x},\vec{y})=\big\{0\leq i\leq n:\text{~there exists $j$ that $0\leq j\leq n$ and~}y_i=x_j\big\}
\]
and
\begin{align*}
\zeta(\vec{x},\vec{y})=\big\{&0\leq i\leq n-1:\text{~there exists $j$ that $0\leq j\leq n-1$ and~}\\
&y_i=x_j\text{~while~}y_{i+1}=x_{j+1}\}.
\end{align*}

Let $\{\alpha_n\}_{n\geq 0}$ be a self-avoiding random walk on $\Gamma_2$ that $\alpha_0\in \Gamma_3$ and
\[
\widetilde{P}\Big(\alpha_n-\alpha_{n-1}=e_j\Big|\alpha_l,0\leq l\leq n-1\Big)=\frac{1}{\lfloor\frac{d}{\log d}\rfloor}
\]
for each $d-\lfloor\frac{d}{\log d}\rfloor+1\leq j\leq d$ and $n$ that $\lfloor\log d\rfloor\mid n$ while
\[
\widetilde{P}\Big(\alpha_n-\alpha_{n-1}=y\Big|\alpha_l,0\leq l\leq n-1\Big)=\frac{1}{|R(\vec{\alpha},n)|}
\]
for each $n$ that $\lfloor\log d\rfloor\nmid n$ and $y\in R(\vec{\alpha},n)$, where
\[
R(\vec{\alpha},n)=\Big\{u:~u-\alpha_{n-1}\in \big\{\pm e_j:~1\leq j\leq d-\lfloor\frac{d}{\log d}\rfloor\big\}\text{~and~}u\neq \alpha_l \text{~for all~}0\leq l\leq n-1\Big\}
\]
while $\widetilde{P}$ is the probability measure of $\{\alpha_n\}_{n\geq 1}$. We use $\vec{\alpha}_n$ to denote the path $(\alpha_0,\ldots,\alpha_n)$, then it is easy to check that $\vec{\alpha}_n\in L_n(\alpha_0)$ for each $n\geq 1$. Note that
\begin{equation}\label{equ 43.1}
|R(\vec{\alpha},n)|\geq 2(d-\lfloor\frac{d}{\log d}\rfloor)-\lfloor \log d\rfloor>0
\end{equation}
for sufficiently large $d$. This is because $\sum_{k=d-\lfloor \frac{d}{\log d}\rfloor+1}^d y(k)=\sum_{k=d-\lfloor \frac{d}{\log d}\rfloor+1}^d\alpha_n(k)$
for $y\in R(\vec{\alpha},n)$ (where $y(k)$ is the $k$th coordinate of $y$) while
\[
\sum_{k=d-\lfloor \frac{d}{\log d}\rfloor+1}^d\alpha_n(k)>\sum_{k=d-\lfloor \frac{d}{\log d}\rfloor+1}^d\alpha_l(k)
\]
for each $l$ that $n-l>\log d$.

It is easy to check that
\begin{equation}\label{equ 43.0}
\sum_{k=d-\lfloor\log d\rfloor+1}^d\alpha_n(k)=1+\big\lfloor\frac{n}{\lfloor\log d\rfloor}\big\rfloor
\end{equation}
for each $n\geq 1$ according to the definition of $\alpha_n$ and $\Gamma_3$. This property will be utilized repeatedly in the proof of Lemma \ref{lemma 4.4}.

Let $\{q_n\}_{n\geq 0}$ be an independent copy of $\{\alpha_n\}_{n\geq 0}$ and $\vec{q}_n=(q_0,\ldots,q_n)$, then we define
\[
\sigma(n)=\sigma(\vec{\alpha}_n,\vec{q}_n)=\big\{0\leq i\leq n:\text{~there exists $j$ that $0\leq j\leq n$ and~}q_i=\alpha_j\big\}
\]
and
\begin{align*}
\zeta(n)=\zeta(\vec{\alpha}_n,\vec{q}_n)=\big\{0\leq i\leq n-1:&\text{~there exists $j$ that $0\leq j\leq n-1$ and~}\\
&q_i=\alpha_j\text{~while~}q_{i+1}=\alpha_{j+1}\}.
\end{align*}
Furthermore, we define
\[
\sigma=\bigcup_{n\geq 1}\sigma(n)=\big\{ i\geq 0:\text{~there exists $j\geq 0$ that~}q_i=\alpha_j\big\}
\]
and
\[
\zeta=\bigcup_{n\geq 1}\zeta(n)=\big\{ i\geq 0:\text{~there exists $j\geq 0$ that~}q_i=\alpha_j\text{~and~}q_{i+1}=\alpha_{j+1}\big\}.
\]
For any $x,y\in \Gamma_3$, we denote by $\widetilde{P}_{x,y}$ the probability measure of $\{\alpha_n,q_n\}_{n\geq 1}$ with $\alpha_0=x$ and $q_0=y$. The expectation operator with respect to $\widetilde{P}_{x,y}$ is denoted by $\widetilde{E}_{x,y}$. The follow lemma is crucial for us to prove Lemma \ref{lemma 4.4}.
\begin{lemma}\label{lemma 43.1}
For any $B\subseteq \Gamma_3$,
\[
P_{\lambda,d}\Big(D_3(n,B)\neq \emptyset\text{~i.o.}\Big)\geq \frac{1}{\frac{1}{|B|^2}\sum_{x\in B}\sum_{y\in B}\widetilde{E}_{x,y}\Big[M_2^{^{|\sigma\setminus\zeta|}}(\frac{2d+\lambda\Theta}{\lambda E\rho})^{|\zeta|}\Big]},
\]
where $M_2=\frac{3\Theta^2}{(E\rho)^2}$.
\end{lemma}

The proof of Lemma \ref{lemma 43.1} is given at the end of this subsection. Now we show how to utilize Lemma \ref{lemma 43.1} to prove Lemma \ref{lemma 4.4}.

\proof[Proof of Lemma \ref{lemma 4.4}]

We define $\kappa=\inf\{i\geq 0:\text{~there exists $j\geq 0$ that~}q_i=\alpha_j\}$. If $\kappa=+\infty$, then $|\sigma|=|\zeta|=0$. As a result,
\begin{align}\label{equ 43.2}
\widetilde{E}_{x,y}\Big[M_2^{^{|\sigma\setminus\zeta|}}(\frac{2d+\lambda\Theta}{\lambda E\rho})^{|\zeta|}\Big]=&\widetilde{P}_{x,y}(\kappa=+\infty)\\
&+\widetilde{E}_{x,y}\Bigg(\widetilde{E}_{x,y}\Big(M_2^{^{|\sigma\setminus\zeta|}}(\frac{2d+\lambda\Theta}{\lambda E\rho})^{|\zeta|}\Big|\kappa<+\infty\Big)1_{\{\kappa<+\infty\}}\Bigg).\notag
\end{align}
We claim that there exists $M_4>0$ which does not depend on $d$ that
\begin{equation}\label{equ 43.5}
\widetilde{P}_{x,y}(\kappa<+\infty)\leq \frac{M_4(\log d)^2}{d}
\end{equation}
for any $x,y\in \Gamma_3, x\neq y$. The proof of Equation \eqref{equ 43.5} will be given later.

Reference \cite{Xue2017b} gives a detailed calculation of the upper bound of the function $f(C_1,C_2)=\widetilde{E}_{x,y}\Big(C_1^{^{|\sigma\setminus\zeta|}}C_2^{|\zeta|}\Big)$ for $x=y=O$. For the general case where $(x,y)\neq (O,O)$, the calculation is still valid after modifying some details. According to a similar analysis with that leads to Lemma 3.4 of \cite{Xue2017b},
for any $C_1,C_2>0$, there exists $M_3>0$ which do not depend on $d,C_1,C_2,x,y$ that
\begin{align}\label{equ 43.3}
&\widetilde{E}_{x,y}\Big(C_1^{^{|\sigma\setminus\zeta|}}C_2^{|\zeta|}\Big|\kappa<+\infty\Big)\leq  \notag\\
&C_1+\frac{C_1C_2}{\lfloor\frac{d}{\log d}\rfloor}\sum_{k=1}^{+\infty}\Big(\frac{C_2}{2(d-\lfloor\frac{d}{\log d}\rfloor)-\lfloor\log d\rfloor}+\frac{C_2}{\lfloor\frac{d}{\log d}\rfloor\lfloor\log d\rfloor^3}+\frac{M_3(\log d)^5C_1}{d}\Big)^k.
\end{align}

Let $C_1=M_2$ and $C_2=\frac{2d+\lambda\Theta}{\lambda E\rho}$ for $\lambda>\frac{1}{E\rho}$, then
\[
\lim_{d\rightarrow+\infty}\frac{C_2}{2(d-\lfloor\frac{d}{\log d}\rfloor)-\lfloor\log d\rfloor}+\frac{C_2}{\lfloor\frac{d}{\log d}\rfloor\lfloor\log d\rfloor^3}+\frac{M_3(\log d)^5C_1}{d}=\frac{1}{\lambda E\rho}<1.
\]
We choose $c\in (\frac{1}{\lambda E\rho},1)$, then for sufficiently large $d$,
\[
\frac{\frac{2d}{\lambda E\rho}}{2(d-\lfloor\frac{d}{\log d}\rfloor)-\lfloor\log d\rfloor}+\frac{\frac{2d}{\lambda E\rho}}{\lfloor\frac{d}{\log d}\rfloor\lfloor\log d\rfloor^3}+\frac{M_3(\log d)^5M_2}{d}\leq c
\]
and
\begin{equation}\label{equ 43.4}
\widetilde{E}_{x,y}\Big(M_2^{^{|\sigma\setminus\zeta|}}(\frac{2d+\lambda\Theta}{\lambda E\rho})^{|\zeta|}\Big|\kappa<+\infty\Big)
\leq M_2\big[1+\frac{4\log d}{(1-c)\lambda E\rho}\big]
\end{equation}
by Equation \eqref{equ 43.3}.

By Equations \eqref{equ 43.2}, \eqref{equ 43.5} and \eqref{equ 43.4}, for sufficiently large $d$ and any $B\subseteq \Gamma_3$,
\begin{equation}\label{equ 43.6}
\widetilde{E}_{x,y}\Big[M_2^{^{|\sigma\setminus\zeta|}}(\frac{2d+\lambda\Theta}{\lambda E\rho})^{|\zeta|}\Big]\leq 1+M_2\big[1+\frac{4\log d}{(1-c)\lambda E\rho}\big] \frac{M_4(\log d)^2}{d}
\end{equation}
for any $x,y\in B,x\neq y$.

By Equations \eqref{equ 43.4} and  \eqref{equ 43.6}, for sufficiently large $d$ and any $B\subseteq \Gamma_3$ with $|B|=\lfloor d^{1/4}\rfloor$,
\begin{align*}
&\frac{1}{|B|^2}\sum_{x\in B}\sum_{y\in B}\widetilde{E}_{x,y}\Big[M_2^{^{|\sigma\setminus\zeta|}}(\frac{2d+\lambda\Theta}{\lambda E\rho})^{|\zeta|}\Big]\\
&=\frac{1}{|B|^2}\Bigg(\sum_{x\in B}\widetilde{E}_{x,x}\Big[M_2^{^{|\sigma\setminus\zeta|}}(\frac{2d+\lambda\Theta}{\lambda E\rho})^{|\zeta|}\Big]
+\sum_{x\in B}\sum_{y:y\neq x}\widetilde{E}_{x,y}\Big[M_2^{^{|\sigma\setminus\zeta|}}(\frac{2d+\lambda\Theta}{\lambda E\rho})^{|\zeta|}\Big]\Bigg)\\
&\leq \frac{|B| M_2\big[1+\frac{4\log d}{(1-c)\lambda E\rho}\big]}{|B|^2}+\frac{|B|(|B|-1)\big(1+M_2\big[1+\frac{4\log d}{(1-c)\lambda E\rho}\big] \frac{M_4(\log d)^2}{d}\big)}{|B|^2}\\
&=\frac{M_2\big[1+\frac{4\log d}{(1-c)\lambda E\rho}\big]}{\lfloor d^{1/4}\rfloor}+\frac{(\lfloor d^{1/4}\rfloor-1)\big(1+M_2\big[1+\frac{4\log d}{(1-c)\lambda E\rho}\big] \frac{M_4(\log d)^2}{d}\big)}{\lfloor d^{1/4}\rfloor}
\end{align*}
and hence
\[
P_{\lambda,d}\Big(D_3(n,B)\neq \emptyset\text{~i.o.}\Big)\geq \frac{1}{\frac{M_2\big[1+\frac{4\log d}{(1-c)\lambda E\rho}\big]}{\lfloor d^{1/4}\rfloor}+\frac{(\lfloor d^{1/4}\rfloor-1)\big(1+M_2\big[1+\frac{4\log d}{(1-c)\lambda E\rho}\big] \frac{M_4(\log d)^2}{d}\big)}{\lfloor d^{1/4}\rfloor}}
\]
by Lemma \ref{lemma 43.1}. Then according to the definition of $\Delta(d)$,
\[
\Delta(d)\geq \frac{1}{\frac{M_2\big[1+\frac{4\log d}{(1-c)\lambda E\rho}\big]}{\lfloor d^{1/4}\rfloor}+\frac{(\lfloor d^{1/4}\rfloor-1)\big(1+M_2\big[1+\frac{4\log d}{(1-c)\lambda E\rho}\big] \frac{M_4(\log d)^2}{d}\big)}{\lfloor d^{1/4}\rfloor}}
\]
and hence
\[
\lim_{d\rightarrow+\infty}\Delta(d)=1.
\]
To finish this proof, we only need to show that Equation \eqref{equ 43.5} holds. According to the definition of $\kappa$,
\[
\{\kappa<\lfloor \log d\rfloor\}=\{\exists~i<\lfloor\log d\rfloor,j\geq0\text{~that~}\alpha_j=q_i\}
\]
while
\[
\{\lfloor\log d\rfloor\leq \kappa<+\infty\}\subseteq \{\exists~i\geq \lfloor\log d\rfloor,j\geq 0\text{~that~}\alpha_j=q_i\}.
\]
Therefore, for $x\neq y$,
\begin{align}\label{equ 43.7}
\widetilde{P}_{x,y}(\kappa<+\infty)\leq &\widetilde{P}_{x,y}(\exists~i<\lfloor\log d\rfloor,j\geq0\text{~that~}\alpha_j=q_i)\\
&+\widetilde{P}_{x,y}(
\exists~i\geq \lfloor\log d\rfloor,j\geq 0\text{~that~}\alpha_j=q_i).\notag
\end{align}

For $x,y\in \Gamma_3, x\neq y$,
\[
\widetilde{P}_{x,y}(\exists~i<\lfloor\log d\rfloor,j\geq0\text{~that~}\alpha_j=q_i)
\leq \sum_{j=1}^{+\infty}\widetilde{P}_{x,y}(\alpha_j=x)+\sum_{i=1}^{\lfloor\log d\rfloor-1}\sum_{j\geq 0}\widetilde{P}_{x,y}(q_i=\alpha_j).
\]
According to the definition of $\Gamma_3$ and Equation \eqref{equ 43.0}, $\sum_{k=d-\lfloor\frac{d}{\log d}\rfloor+1}^dx(k)=1$ while
\[\sum_{k=d-\lfloor\frac{d}{\log d}\rfloor+1}^d\alpha_j(k)\geq 2\] for $j\geq \lfloor\log d\rfloor$. Hence, by Equation \eqref{equ 43.1},
\begin{align}\label{equ 43.8}
\sum_{j=1}^{+\infty}\widetilde{P}_{x,y}(\alpha_j=x)&=\sum_{j=1}^{\lfloor\log d\rfloor-1}\widetilde{P}_{x,y}(\alpha_j=x)\leq \sum_{j=1}^{\lfloor\log d\rfloor-1}\sup_{u\in \mathbb{Z}^d}\widetilde{P}_{x,y}(\alpha_j=u)\\
&=\sum_{j=1}^{\lfloor\log d\rfloor-1}\widetilde{E}_{x,y}\big(\frac{1}{|R(\vec{\alpha},n)|}\big)\leq \frac{\lfloor\log d\rfloor-1}{2(d-\lfloor\frac{d}{\log d}\rfloor)-\lfloor \log d\rfloor}.\notag
\end{align}
According to Equation \eqref{equ 43.0}, $\sum_{k=d-\lfloor\frac{d}{\log d}\rfloor+1}^d\alpha_j(k)\neq \sum_{k=d-\lfloor\frac{d}{\log d}\rfloor+1}^dq_i(k)$
when $|i-j|>\lfloor\log d\rfloor$, therefore by Equation \eqref{equ 43.1},
\begin{align}\label{equ 43.9}
&\sum_{i=1}^{\lfloor\log d\rfloor-1}\sum_{j\geq 0}\widetilde{P}_{x,y}(q_i=\alpha_j)\notag\\
&=\sum_{i=1}^{\lfloor\log d\rfloor-1} \sum_{j:|j-i|\leq \lfloor \log d\rfloor}\widetilde{P}_{x,y}(q_i=\alpha_j)\leq \sum_{i=1}^{\lfloor\log d\rfloor-1}\sum_{j:|j-i|\leq \lfloor \log d\rfloor}\sup_{u\in \mathbb{Z}^d}\widetilde{P}_{x,y}(q_i=u) \\
&=\sum_{i=1}^{\lfloor\log d\rfloor-1}\sum_{j:|j-i|\leq \lfloor \log d\rfloor}\widetilde{E}_{x,y}\big(\frac{1}{|R(\vec{q},n)|}\big)\leq \frac{(\lfloor\log d\rfloor-1)(2\lfloor\log d\rfloor+1)}{2(d-\lfloor\frac{d}{\log d}\rfloor)-\lfloor \log d\rfloor}.\notag
\end{align}
By Equations \eqref{equ 43.8} and \eqref{equ 43.9},
\begin{equation}\label{equ 43.10}
\widetilde{P}_{x,y}(\exists~i<\lfloor\log d\rfloor,j\geq0\text{~that~}\alpha_j=q_i)\leq \frac{2(\log d)^2}{d}
\end{equation}
for sufficiently large $d$.

For $w=(w_1,\ldots,w_d)\in \mathbb{Z}^d$, we use $\vartheta(w)$ to denote
\[
(w_{d-\lfloor\frac{d}{\log d}\rfloor+1},\ldots,w_d)\in \mathbb{Z}^{\lfloor\frac{d}{\log d}\rfloor}.
\]
Then, by Equation \eqref{equ 43.0},
\begin{equation}\label{equ 43.11}
\{\exists~i\geq \lfloor\log d\rfloor,j\geq 0\text{~that~}\alpha_j=q_i\}\subseteq \big\{\vartheta(\alpha_{_{k\lfloor\log d\rfloor}})=\vartheta(q_{_{k\lfloor\log d\rfloor}})\text{~for some~}k\geq 1\big\}.
\end{equation}
According to the definitions of $\alpha_n$ and $q_n$, $\{\vartheta(\alpha_{k\lfloor\log d\rfloor})\}_{k\geq 1}$ and $\{\vartheta(q_{k\lfloor\log d\rfloor})\}_{k\geq 1}$ are two independent oriented random walks on $\mathbb{Z}^{\lfloor\frac{d}{\log d}\rfloor}$. Then, according to the lemma given in \cite{Cox1983} about the first collision time of two independent oriented random walks on the lattice, there exists $M_5$ which does not depend on $d,x,y$ that
\begin{equation}\label{equ 43.12}
\widetilde{P}_{x,y}\big(\vartheta(\alpha_{_{k\lfloor\log d\rfloor}})=\vartheta(q_{_{k\lfloor\log d\rfloor}})\text{~for some~}k\geq 1\big)\leq \frac{M_5}{\lfloor\frac{d}{\log d}\rfloor}.
\end{equation}
By Equations \eqref{equ 43.11} and \eqref{equ 43.12},
\begin{equation}\label{equ 43.13}
\widetilde{P}_{x,y}\big(\exists~i\geq \lfloor\log d\rfloor,j\geq 0\text{~that~}\alpha_j=q_i\big)\leq \frac{M_5}{\lfloor\frac{d}{\log d}\rfloor}.
\end{equation}
By equations \eqref{equ 43.7}, \eqref{equ 43.10} and \eqref{equ 43.13},
\begin{equation}\label{equ 43.14}
\widetilde{P}_{x,y}(\kappa<+\infty)\leq \frac{2(\log d)^2}{d}+\frac{M_5}{\lfloor\frac{d}{\log d}\rfloor}
\end{equation}
for sufficiently large $d$, where $M_5$ does not depend on $d,x,y$. Equation \eqref{equ 43.5} follows from Equation \eqref{equ 43.14} directly.

\qed

At the end of this subsection, we give the proof of Lemma \ref{lemma 43.1}. The proof utilizes the following Proposition given in \cite{Xue2017b}.
\begin{proposition}\label{proposition 43.1}
If $A_1, A_2,\ldots, A_n$ are $n$ arbitrary random events defined under the same probability space such that $P(A_i)>0$ for $1\leq i\leq n$ and $q_1,q_2,\ldots,q_n$ are $n$ positive constants such that
$\sum_{j=1}^nq_j=1$, then
\[
P(\bigcup_{j=1}^{+\infty}A_j)\geq \frac{1}{\sum\limits_{i=1}^n\sum\limits_{j=1}^nq_iq_j\frac{P(A_i\bigcap A_j)}{P(A_i)P(A_j)}}.
\]
\end{proposition}
This proposition is Lemma 3.3 of \cite{Xue2017b} and a detailed proof is given there.

\proof[Proof of Lemma \ref{lemma 43.1}]

For any $B\subseteq \Gamma_3$ and $\vec{l}\in L_n(B)$, we denote by $G_{\vec{l}}$ the event that $\vec{l}$ is an infectious path. According to the definition of $D_3(n,B)$, it is easy to check that
\[
\big\{D_3(n,B)\neq \emptyset\text{~i.o.}\big\}\supseteq \bigcap_{n=1}\bigcup_{\vec{l}\in L_n(B)}G_{\vec{l}}.
\]
As a result,
\begin{equation}\label{equ 43.15}
P_{\lambda,d}\big(D_3(n,B)\neq \emptyset\text{~i.o.}\big)\geq \lim_{n\rightarrow+\infty}P_{\lambda,d}\big(\bigcup_{\vec{l}\in L_n(B)}G_{\vec{l}}\big),
\end{equation}
since $\bigcup_{\vec{l}\in L_n(B)}G_{\vec{l}}\subseteq \bigcup_{\vec{l}\in L_m(B)}G_{\vec{l}}$ for $n>m$.

For each $\vec{l}=(l_0,\ldots,l_n)\in L_n(B)$, we define \[g_{\vec{l}}=\frac{\widetilde{P}_{l_0}\big(\vec{\alpha}_n=\vec{l}\big)}{|B|},\] where
\[
\vec{\alpha}_n=(\alpha_0,\ldots,\alpha_n),
\]
then it easy to check that $\sum_{\vec{l}\in B}g_{\vec{l}}=1$. Then, by Proposition \ref{proposition 43.1},
\begin{align}\label{equ 43.16}
P_{\lambda,d}\big(\bigcup_{\vec{l}\in L_n(B)}G_{\vec{l}}\big)&\geq \frac{1}{\sum_{\vec{x}\in L_n(B)}\sum_{\vec{y}\in L_n(B)}g_{\vec{x}}g_{\vec{y}}\frac{P_{\lambda,d}(G_{\vec{x}}\bigcap G_{\vec{y}})}{P_{\lambda,d}(G_{\vec{x}})P_{\lambda,d}(G_{\vec{y}})}} \notag\\
&=\frac{1}{\frac{1}{|B|^2}\sum\limits_{x\in B}\sum\limits_{y\in B}\sum\limits_{\vec{x}\in L_n(x)}\sum\limits_{\vec{y}\in L_n(y)}\frac{\widetilde{P}_{x}\big(\vec{\alpha}_n=\vec{x}\big)\widetilde{P}_{y}\big(\vec{q}_n=\vec{y}\big)P_{\lambda,d}(G_{\vec{x}}\bigcap G_{\vec{y}})}{P_{\lambda,d}(G_{\vec{x}})P_{\lambda,d}(G_{\vec{y}})}}.
\end{align}

 Now we deal with the factor $\frac{P_{\lambda,d}(G_{\vec{x}}\bigcap G_{\vec{y}})}{P_{\lambda,d}(G_{\vec{x}})P_{\lambda,d}(G_{\vec{y}})}$. According to our assumption of the model, the denominator
 \[
 P_{\lambda,d}(G_{\vec{x}})P_{\lambda,d}(G_{\vec{y}})=\prod_{l=0}^{n-1}P_{\lambda,d}\big(U(x_l,x_{l+1})<Y(x_l)\big)\prod_{l=1}^{n-1} P_{\lambda,d}\big(U(y_l,y_{l+1})<Y(y_l)\big),
 \]
since both $\vec{x}$ and $\vec{y}$ are self-avoiding. According to our assumption of the model, for any $y_i\not \in \vec{x}$ and $x_j\not\in \vec{y}$, the numerator $P_{\lambda,d}(G_{\vec{x}}\bigcap G_{\vec{y}})$ has factors $P_{\lambda,d}\big(U(x_j,x_{j+1})<Y(x_j)\big)$ and $P_{\lambda,d}\big(U(y_i,y_{i+1})<Y(y_j)\big)$, which can be cancelled with the same factors in the denominator. For each $l\in \zeta(\vec{x},\vec{y})$, there exists $0\leq k\leq n-1$ and $u,v\in\mathbb{Z}^d$ that $y_l=x_k=u$ and $y_{l+1}=x_{k+1}=v$. Then, the numerator $P_{\lambda,d}(G_{\vec{x}}\bigcap G_{\vec{y}})$ has the factor $P_{\lambda,d}\big(U(u,v)<Y(u)\big)$ while the denominator $P_{\lambda,d}(G_{\vec{x}})P_{\lambda,d}(G_{\vec{y}})$ has the factor $\Big[P_{\lambda,d}\big(U(u,v)<Y(u)\big)\Big]^2$. Hence, $\frac{P_{\lambda,d}(G_{\vec{x}}\bigcap G_{\vec{y}})}{P_{\lambda,d}(G_{\vec{x}})P_{\lambda,d}(G_{\vec{y}})}$ has the factor
\[
\frac{P_{\lambda,d}\big(U(u,v)<Y(u)\big)}{\Big[P_{\lambda,d}\big(U(u,v)<Y(u)\big)\Big]^2}=\frac{2d}{\lambda E\big(\frac{\rho}{1+\frac{\lambda\rho}{2d}}\big)}
\]
for each $l\in \zeta(\vec{x},\vec{y})$.

For each $m\in \sigma(\vec{x},\vec{y})\setminus \zeta(\vec{x},\vec{y})$, there exists $0\leq r\leq n-1$ and $u_1,v_1,v_2\in \mathbb{Z}^d, v_1\neq v_2$ that $y_m=x_r=u_1$, $y_{m+1}=v_2$ and $x_{r+1}=v_1$. Then, the denominator $P_{\lambda,d}(G_{\vec{x}})P_{\lambda,d}(G_{\vec{y}})$ has the factor $P_{\lambda,d}\big(U(u_1,v_1)<Y(u_1)\big)P_{\lambda,d}\big(U(u_1,v_2)<Y(u_1)\big)$ while the numerator has a factor at most $P_{\lambda,d}\big(\widetilde{U}(u_1,v_1)<Y(u_1),\widetilde{U}(u_1,v_2)<Y(u_1)\big)$, where $\widetilde{U}(u_1,v_1),\widetilde{U}(u_1,v_2)$ are independent exponential times with rate $\frac{\lambda \Theta}{2d}$ and are independent with $Y(u_1)$. Note that here we replace $U(u_1,v_2)$ and $U(u_1,v_1)$ by $\widetilde{U}(u_1,v_2)$ and $\widetilde{U}(u_1,v_1)$ in case some other event in the numerator depends on the exponential time $U(v_1,u_1)$ or $U(v_2,u_1)$, which are independent with $U(u_1,v_2)$ and $U(u_1,v_1)$ under the quenched measure but positively correlated under the annealed measure. As a result, for each $m\in \sigma(\vec{x},\vec{y})$, $\frac{P_{\lambda,d}(G_{\vec{x}}\bigcap G_{\vec{y}})}{P_{\lambda,d}(G_{\vec{x}})P_{\lambda,d}(G_{\vec{y}})}$ has the factor at most
\[
\frac{P_{\lambda,d}\big(\widetilde{U}(u_1,v_1)<Y(u_1),\widetilde{U}(u_1,v_2)<Y(u_1)\big)}{P_{\lambda,d}\big(U(u_1,v_1)<Y(u_1)\big)P_{\lambda,d}\big(U(u_1,v_2)<Y(u_1)\big)} =\frac{2\Theta^2}{(1+\frac{\lambda\Theta}{2d})(1+\frac{\lambda\Theta}{d})\big[E(\frac{\rho}{1+\frac{\lambda\rho}{2d}})\big]^2}.
\]
Inclusion,
\begin{align}\label{equ 43.17}
&\frac{P_{\lambda,d}(G_{\vec{x}}\bigcap G_{\vec{y}})}{P_{\lambda,d}(G_{\vec{x}})P_{\lambda,d}(G_{\vec{y}})}\notag\\
&\leq\Big[\frac{2d}{\lambda E\big(\frac{\rho}{1+\frac{\lambda\rho}{2d}}\big)}\Big]^{|\zeta(\vec{x},\vec{y})|}\Big[\frac{2\Theta^2}{(1+\frac{\lambda\Theta}{2d})(1+\frac{\lambda\Theta}{d})
\big[E(\frac{\rho}{1+\frac{\lambda\rho}{2d}})\big]^2}\Big]^{|\sigma(\vec{x},\vec{y})\setminus\zeta(\vec{x},\vec{y})|}\\
&\leq \Big[\frac{2d+\lambda\Theta}{\lambda E\rho}\Big]^{|\zeta(\vec{x},\vec{y})|}\Big[\frac{3\Theta^2}{(E\rho)^2}\Big]^{|\sigma(\vec{x},\vec{y})\setminus\zeta(\vec{x},\vec{y})|} \notag
\end{align}
for sufficiently large $d$.

By Equations \eqref{equ 43.16}, \eqref{equ 43.17} and the definitions of $\sigma(n), \zeta(n)$,
\begin{align*}
&P_{\lambda,d}\big(\bigcup_{\vec{l}\in L_n(B)}G_{\vec{l}}\big)\geq  \\
&\frac{1}{\frac{1}{|B|^2}\sum\limits_{x\in B}\sum\limits_{y\in B}\sum\limits_{\vec{x}\in L_n(x)}\sum\limits_{\vec{y}\in L_n(y)}\widetilde{P}_{x}\big(\vec{\alpha}_n=\vec{x}\big)\widetilde{P}_{y}\big(\vec{q}_n=\vec{y}\big)\Big[\frac{2d+\lambda\Theta}{\lambda E\rho}\Big]^{|\zeta(\vec{x},\vec{y})|}\Big[\frac{3\Theta^2}{(E\rho)^2}\Big]^{|\sigma(\vec{x},\vec{y})\setminus\zeta(\vec{x},\vec{y})|}}\\
&=\frac{1}{\frac{1}{|B|^2}\sum\limits_{x\in B}\sum\limits_{y\in B}\widetilde{E}_{x,y}\Big(\big(\frac{2d+\lambda\Theta}{\lambda E\rho}\big)^{|\zeta(n)|}\big(\frac{3\Theta^2}{(E\rho)^2}\big)^{|\sigma(n)\setminus\zeta(n)|}\Big)}.
\end{align*}
Then, according to the definition of $\sigma$ and $\zeta$,
\begin{equation}\label{equ 43.18}
\lim_{d\rightarrow+\infty}P_{\lambda,d}\big(\bigcup_{\vec{l}\in L_n(B)}G_{\vec{l}}\big)\geq \frac{1}{\frac{1}{|B|^2}\sum\limits_{x\in B}\sum\limits_{y\in B}\widetilde{E}_{x,y}\Big(\big(\frac{2d+\lambda\Theta}{\lambda E\rho}\big)^{|\zeta|}\big(\frac{3\Theta^2}{(E\rho)^2}\big)^{|\sigma\setminus\zeta|}\Big)}.
\end{equation}

Lemma \ref{lemma 43.1} follows directly from Equations \eqref{equ 43.15} and \eqref{equ 43.18}. 

\qed

\section{Proof of Theorem \ref{theorem 2.2}}\label{section 5}
In this section, we give the proof of Theorem \ref{theorem 2.2}.

\proof[Proof of Theorem \ref{theorem 2.2}]
For given $\lambda>\frac{1}{E\rho}$, according to Theorem \ref{theorem 2.1 main},
\[
P_{\lambda,d}\big(C_t^O\neq \emptyset,\forall~t\geq 0\big)>0
\]
for sufficiently large $d$. Therefore, $\lambda_c(d)\leq \lambda$ for sufficiently large $d$ and hence
\[
\limsup_{d\rightarrow+\infty}\lambda_c(d)\leq \lambda.
\]
Let $\lambda\rightarrow\frac{1}{E\rho}$, we have
\begin{equation}\label{equ 5.1}
\limsup_{d\rightarrow+\infty}\lambda_c(d)\leq \frac{1}{E\rho}.
\end{equation}
According to a similar analysis with that leads to Equation \eqref{equ 5.1},
\begin{equation}\label{equ 5.2}
\limsup_{d\rightarrow+\infty}\beta_c(d)\leq \frac{1}{E\rho}.
\end{equation}

Now we let $\lambda=\frac{\gamma}{E\rho}$ for given $\gamma<1$. For a self-avoiding path
\[\vec{l}=(O,l_1,\ldots,l_n)\] on $\mathbb{Z}^d$ starting at $O$ with length $n$,
\begin{align*}
&P_{\lambda,d}\big(\vec{l}\text{~is an infectious path}\big)=E_{\mu_d}\Big[P_{\lambda,\omega}\big(\vec{l}\text{~is an infectious path}\big)\Big]\\
&=E_{\mu_d}\Big[\prod_{i=0}^{n-1}P_{\lambda,\omega}\big(U(l_i,l_{i+1})<Y(l_i)\big)\Big]
=E_{\mu_d}\Big[\prod_{i=0}^{n-1}\frac{\frac{\gamma\rho(l_i,l_{i+1},\omega)}{2dE\rho}}{\frac{\gamma\rho(l_i,l_{i+1},\omega)}{2dE\rho}+1}\Big]\\
&=\Big[\frac{\gamma}{2dE\rho}E\big(\frac{\rho}{1+\frac{\gamma\rho}{2dE\rho}}\big)\Big]^n\leq \frac{\gamma^n}{(2d)^n}.
\end{align*}
For a self-avoiding path on $\mathbb{Z}^d$, each step has at most $2d$ choices. Therefore, the number of self-avoiding paths with length $n$ staring at $O$ is at most $(2d)^n$. As a result,
\begin{align*}
&P_{\lambda,n}\big(\text{~there exists an infectious path with length $n$ starting at $O$}\big)\\
&\leq (2d)^n\frac{\gamma^n}{(2d)^n}=\gamma^n.
\end{align*}
Since $\sum_{n=0}^{+\infty}\gamma^n<+\infty$ for $\gamma<1$, according to the Borel-Cantelli's lemma,
\begin{equation}\label{equ 5.3}
P_{\lambda,n}\big(\text{~there exist arbitrary long infectious paths starting at $O$}\big)=0.
\end{equation}
By Lemma \ref{lemma 4.1},
\begin{align*}
&\{I_t^O\neq \emptyset,\forall~t\geq 0\}\\
&\subseteq \{\text{~there exist arbitrary long infectious paths starting at $O$}\},
\end{align*}
since infectious vertices never die out when and only when there are infinite many vertices have ever been infected. Then by Equation \eqref{equ 5.3},
\[
P_{\lambda,n}\big(I_t^O\neq \emptyset,\forall~t\geq 0\big)=0
\]
for $\lambda=\frac{\gamma}{E\rho}$ with $\gamma<1$ and hence
\[
\beta_c(d)\geq \frac{\gamma}{E\rho}
\]
for any $\gamma<1$. Let $\gamma\rightarrow 1$, we have $\beta_c(d)\geq \frac{1}{E\rho}$ and hence
\begin{equation}\label{equ 5.4}
\lim_{d\rightarrow+\infty}\beta_c(d)=\frac{1}{E\rho}
\end{equation}
combined with Equation \eqref{equ 5.2}. Theorem \ref{theorem 2.2} follows from Equations \eqref{equ 5.1} and \eqref{equ 5.4} directly.

\qed

\quad

\textbf{Acknowledgments.}
The author is grateful to the financial
support from the National Natural Science Foundation of China with
grant number 11501542 and the financial support from Beijing
Jiaotong University with grant number KSRC16006536.

{}

\begin{thebibliography}{}
\bibitem{Ber2011}Bertacchi, D., Lanchier, N. and Zucca, F. (2011). Contact and voter processes on the infinite percolation cluster as models of host-symbiont interactions. \emph{The Annals of Applied Probability} \textbf{21}, 1215-1252.
\bibitem{Chen2009}Chen, XX. and Yao, Q. (2009). The complete convergence theorem holds for contact processes on open clusters of $Z^d\times Z^+$. \emph{Journal of Statistical Physics} \textbf{135}, 651-680.
\bibitem{Cox1983}Cox, J. T. and Durrett, R. (1983). Oriented percolation in dimensions $d\geq 4$: bounds and asymptotic formulas.
\emph{Mathematical Proceedings of the Cambridge Philosophical
Society} \textbf{93}, 151-162.
\bibitem{Grif1983}Griffeath, D. (1983). The Binary Contact Path Process. \emph{The Annals of Probability} \textbf{11} 692-705.
\bibitem{Har1974}Harris, T. E. (1974). Contact interactions on a lattice. \emph{The Annals of Probability} \textbf{2}, 969-988.
\bibitem{Hol1981}Holley, R. and Liggett, T. M. (1981). Generalized potlatch and smoothing processes.
\emph{Zeitschrift f\"{u}r Wahrscheinlichkeitstheorie und Verwandte Gebiete } \textbf{55}, 165-195.
\bibitem{Kesten1991}Kesten, H. (1991). Asymptotics in high dimension for the Fortuin-Kasteleyn cluster model. In \emph{Spatial Stochastic Processes}, 57-85. Birkh\"{a}user, Boston.
\bibitem{Lig1985}Liggett, T. M. (1985). \emph{Interacting Particle Systems.} Springer, New York.
\bibitem{Lig1999}Liggett, T. M. (1999). \emph{Stochastic interacting systems: contact, voter and exclusion processes.}
Springer, New York.
\bibitem{Pet2011}Peterson, J. (2011). The contact process on the complete graph with random vertex-dependent infection rates. \emph{Stochastic Processes and their Applications} \textbf{121}(3), 609-629.
\bibitem{Xue2015}Xue, XF. (2015). Contact processes with random vertex weights on oriented lattices. \emph{ALEA-Latin American Journal of Probability and Mathematical Statistics} \textbf{12}, 245-259.
\bibitem{Xue2016}Xue, XF. (2016). Phase transition for the large-dimensional contact process with random recovery rates on open clusters. \emph{Journal of Statistical Physics} \textbf{165}, 845-865.
\bibitem{Xue2017}Xue, XF. (2017). Mean field limit for survival probability of the high-dimensional contact process. \emph{Statistics $\&$ Probability Letters}
\textbf{127}, 178-184.
\bibitem{Xue2017b}Xue, XF. (2017). Asymptotic for critical value of the large-dimensional SIR epidemic on clusters. Arxiv: 1610.09616.
\bibitem{Yao2012}Yao, Q. and Chen, XX. (2012). The complete convergence theorem holds for contact processes in a random environment on $Z^d\times Z^+$. \emph{Stochastic Processes and their Applications} \textbf{122}, 3066-3100.
\end{thebibliography}
\end{document}